\documentclass[11pt, a4paper, fleqn]{article}
\usepackage[latin1]{inputenc}
\usepackage{amsmath}
\usepackage{amsthm}
\usepackage{amsfonts}
\usepackage{amssymb}
\usepackage[T1]{fontenc}
\usepackage{float}
\usepackage[all]{xy}
\usepackage{bbm}
\usepackage{color}
\usepackage{hyperref}
\usepackage{graphicx}
\usepackage{transparent}
\usepackage{graphics}
\usepackage{xcolor}
\usepackage{color}
\usepackage{geometry}
\usepackage{marginnote}

\usepackage{tikz}

\usetikzlibrary{arrows,decorations.pathmorphing,backgrounds,positioning,fit,petri}

\title{$\tau$-Tilting Theory and $\tau$-Slices}
\author{Hipolito Treffinger}

\theoremstyle{plain} 
\newtheorem{theorem}{Theorem}[section]
\newtheorem{prop}[theorem]{Proposition}
\newtheorem{lem}[theorem]{Lemma}
\newtheorem{cor}[theorem]{Corollary}

\theoremstyle{remark}
\newtheorem{rmk}[theorem]{Remark}
\newtheorem{ex}[theorem]{Example}

\theoremstyle{definition}
\newtheorem{defi}[theorem]{Definition}

\newtheorem{question}[theorem]{Question}

\newcommand{\rep}[1]{%
  {%
    \tiny%
    \begin{matrix}%
      #1%
    \end{matrix}%
  }%
}

\usepackage[nottoc]{tocbibind}

\begin{document}

\maketitle

\begin{abstract}
Comparing the module categories of an algebra and of the endomorphism algebra of a given support $\tau$-tilting module, we give a generalization of the Brenner-Butler's tilting theorem in the framework of $\tau$-tilting theory. Afterwards we define $\tau$-slices and prove that complete slices of tilted algebras and local slices of cluster tilted algebras are examples of complete $\tau$-slices. Then we apply this concept to the study of simply connected tilted algebras. Finally, we study the one-point extensions and the split-by-nilpotent extensions of an algebra with $\tau$-slices. 
\end{abstract}

\section{Introduction}

Let $A$ be a finite dimensional basic algebra over an algebraically closed field $k$ and mod$A$ be the category of finitely generated right $A$-modules. The Auslander-Reiten translation in mod$A$ is denoted by $\tau_A$.

Since the end of the twentieth century, tilting theory has played a central r\^ole in the representation theory of finite dimensional algebras. It started with the study of the Coxeter functors defined by Bernstein, Gelfand and Ponomarev in \cite{BGP}, further generalized by Auslander, Platzeck and Reiten in \cite{APR}. Afterwards Brenner and Butler introduced in \cite{BB} tilting functors and proved the so-called tilting theorem. A couple of years later, Happel and Ringel defined tilted algebras as endomorphism algebras of the tilting modules over hereditary algebras. They showed in \cite{HR} that many of the characteristics of the module category of a tilted algebra can be deduced from those of a corresponding hereditary algebra via the tilting theorem. Moreover they showed that given a tilted algebra one can recover the corresponding hereditary algebra as endomorphism algebra of a particular tilting module. They called this tilting module a complete slice.

In the early 2000, Fomin and Zelevinski defined in \cite{FZ1} cluster algebras. These algebras are connected to a wide variety of areas in mathematics, including representation theory. The definition of cluster algebras relies strongly onto the concept of mutation. In \cite{AIR}, Adachi, Iyama and Reiten applied the mutation process to representation theory of algebras and extended classical tilting theory via the so-called $\tau$-tilting theory.

Tilting theory is based on comparing the representation theory of an algebra with that of the endomorphism algebra of a tilting module over that algebra. It is shown in \cite{AIR} that tilting modules are exactly the faithful support $\tau$-tilting modules (see Proposition \ref{lista}). Hence, extending results from tilting theory to $\tau$-tilting theory increases the amount of module categories that can be compared to a given one. The first part of this paper is devoted to the comparison of the module category of a given algebra with the module category of the endomorphism algebra of a given support $\tau$-tilting module. We denote by $FacM$ (or $SubM$) the full subcategory of mod$A$ having as objects the quotient modules (or submodules respectively) of finite direct sums of $M$. Our first result is the following.

\begin{theorem}[Theorem \ref{cociente}]
Let $A$ be an algebra and $M$ a support $\tau$-tilting $A$-module, $B=End_A M$ be its endomorphism algebra and $C=A/Ann M$.
\begin{enumerate}
\item The algebra morphism $\varphi: A \longrightarrow End_B(_BM)$ defined by $\varphi(a)(m)=ma$ for every $m\in M$ and $a\in A$ induces an isomorphism $C\cong End_B(_BM)$. 
\item There exists a torsion pair $(\mathcal{X},\mathcal{Y})$ in mod$B$ such that the functor $Ext^1_A(M,-): Sub( \tau_A M)\longrightarrow \mathcal{X}$ is an equivalence of categories with quasi-inverse $Tor_1^B(-,M): \mathcal{X}\longrightarrow Sub(\tau_AM)$ if and only if $\tau_C M=\tau_A M$.
\end{enumerate}
\end{theorem}

We combine the previous result with one proved by Jasso (\cite[Proposition 3.5]{J}) in order to generalize the tilting theorem of Brenner and Butler. Our statement is the following.

\begin{cor}[Corollary \ref{tauBB}]
Let $A$ be an algebra and $M$ a support $\tau$-tilting $A$-module. Consider the pair of subcategories $(Fac M, Sub(\tau_A M))$ of mod$A$. Let $B=End_A(M)$ be the endomorphism algebra of $M$ and $C= A/Ann M$. Then there exists a torsion pair $(\mathcal{X},\mathcal{Y})$ in mod$B$ such that
\begin{enumerate}
\item The functor $Hom_A(M,-): Fac M\longrightarrow \mathcal{Y}$ is an equivalence of categories with quasi-inverse $-\otimes_B M: \mathcal{Y}\longrightarrow Fac M$.
\end{enumerate}
Moreover $\tau_C M=\tau_A M$ if and only if
\begin{enumerate}
\item[2.] The functor $Ext^1_A(M,-): Sub( \tau_A M)\longrightarrow \mathcal{X}$ is an equivalence of categories with quasi-inverse $Tor_1^B(-,M): \mathcal{X}\longrightarrow Sub(\tau_AM)$.
\end{enumerate}
\end{cor}

Next we turn our attention to a special class of support $\tau$-tilting modules, which we call $\tau$-slices. 

\begin{defi}[Definition \ref{tauslice}]
Let $\Sigma$ be a presection (see Definition \ref{presection}) in the Auslander-Reiten quiver $\Gamma_A$ of $A$. Then $\Sigma$ is a \textit{$\tau$-slice} if $M=\bigoplus\limits_{U\in\Sigma} U$ is a support $\tau$-tilting module. It is a \textit{complete $\tau$-slice} if $M$ is a $\tau$-tilting module. 
\end{defi}

The main motivation to introduce $\tau$-slices is that they respect the assumptions of the previous corollary. The $\tau$-slices generalize the complete slices defined in \cite{HR} and the local slices defined in \cite{ABS2}. For instance the endomorphism algebra of a $\tau$-slice is always hereditary (Proposition \ref{endoslice}). Also we are able to prove the following equivalences.

\begin{theorem}[Proposition \ref{proposicioncompleta}, Proposition \ref{clustertilted}]
Let $A$ an algebra. 
\begin{enumerate}
 \item If $A$ is a tilted algebra, then $\Sigma$ is a complete slice if and only if $\Sigma$ is a complete $\tau$-slice.
 \item If $A$ is a cluster-tilted algebra, then $\Sigma$ is a local slice if and only if $\Sigma$ is a complete $\tau$-slice.
\end{enumerate}
\end{theorem}

Moreover, we use $\tau$-slices to get a result on simply connected tilted algebras. 

\begin{theorem}[Theorem \ref{simplyconnectedtilted}]
Let $A$ be an algebra and $\Gamma_A$ be its Auslander-Reiten quiver. Then $\Gamma_A$ has a connected component $\Gamma$ that is simply connected, convex, generalized standard and has a complete $\tau$-slice $\Sigma$ if and only if $A$ is a simply connected tilted algebra having $\Sigma$ as a complete slice in mod$A$.
\end{theorem}

We also show that some one-point extensions (Proposition \ref{onepoint}, Corollary \ref{corex}) and split-by-nilpotent extensions (Theorem \ref{splitex}) preserve $\tau$-slices. 

\subsection{Preliminaries}

In this paper, by algebra we mean a finite dimensional basic algebra over an algebraically closed field $k$. Given an algebra $A$, we consider the category mod$A$ of finitely generated right $A$-modules. For a given object $M$ of mod$A$ we denote by $add M$ (respectively $Fac M$, $Sub M$) the category of all direct summands (respectively quotient modules, submodules) of finite direct sums of copies of $M$. Following the standard notation we let $D=Hom_k(-,k)$ be the duality between mod$A$ and mod$A^{op}$. Moreover the Auslander-Reiten translate of $M$ in mod$A$ is denoted by $\tau_AM$ and the Auslander-Reiten quiver of $A$ is denoted by $\Gamma_A$. If the context is unambiguous, the index of the Auslander-Reiten translation will be omitted leading to the notation $\tau M$. Given a module $X$, the \textit{$\tau$-orbit of $X$} is the class of modules $Y$ such that $Y\cong\tau^nX$for some integer $n\in\mathbb{Z}$. Also, the number of non-isomorphic indecomposable direct summands of $X$ is denoted by $|X|$. We refer the reader to \cite{AsSS} or \cite{ARS} for further details.

Given a path algebra $kQ/I$, the composition of the arrows is denoted by order of appearance, for instance the path $x\overset{\alpha}\rightarrow y\overset{\beta}\rightarrow z$ is written $\alpha\beta$.

Given two modules $X,Y$ in mod$A$ we say that there exists a \textit{path} in mod$A$ from $X$ to $Y$ if there exists a sequence
$$X\overset{f_0}\longrightarrow X_1\overset{f_1}\longrightarrow\dots\overset{f_{n-1}}\longrightarrow X_n\overset{f_n}\longrightarrow Y$$
where the $X_i$ are indecomposable $A$-modules and the $f_i$ are nonzero morphisms.
We then say that $X$ is a \textit{predecessor} of $Y$ and that $Y$ is a \textit{successor} of $X$. Moreover, if the path consists of only one irreducible morphism we say that $X$ is an \textit{immediate predecessor} of $Y$ and that $Y$ is an \textit{immediate successor} of $X$. 

A connected component $\Gamma$ is said to be \textit{generalized standard} if for every pair of modules $X,Y\in \Gamma$ we have rad$^{\infty}_A(X,Y)=0$.

Let $\Sigma$ be a full subquiver of $\Gamma$. Then $\Sigma$ is said to be \textit{convex} if for any $X,Y\in\Sigma$ and every path 
$$X\overset{f_0}\longrightarrow X_1\overset{f_1}\longrightarrow\dots\overset{f_{n-1}}\longrightarrow X_n\overset{f_n}\longrightarrow Y$$
the modules $X_i$ belong to $\Sigma$ for all $i\in\{1,...,n\}$. 

Let $\Gamma$ be a connected component of $\Gamma_A$. A connected full subquiver $\Sigma$ of $\Gamma$ is a \textit{section} if it is acyclic, convex and intersects every $\tau$-orbit of $\Gamma$ exactly once.

The concept of complete slice plays a central part in the present paper. There are many equivalent definitions in the literature. We use the definition given by Ringel in \cite{Rin} for our proofs, which is the following. 

\begin{defi}[\cite{Rin}, Definition 4.2.2]
A finite subquiver $\Sigma\subset\Gamma_A$ is said to be a \textit{slice} if satisfies the following axioms:
\begin{enumerate}
\item $\Sigma$ is a sincere $A$-module.
\item $\Sigma$ is a convex in mod$A$.
\item If $0\longrightarrow L\longrightarrow M\longrightarrow N\longrightarrow 0$ is an almost split sequence in mod $A$, then at most one of $L$ and $N$  belongs to $\Sigma$.
\end{enumerate}
Moreover, $\Sigma$ is called a \textit{complete slice} if it also satisfies the following condition:
\begin{enumerate}
\item[4.] If $0\longrightarrow L\longrightarrow M\longrightarrow N\longrightarrow 0$ is an almost split sequence in mod $A$ such that an indecomposable direct summand $X$ of $M$ is also a direct summand of $\Sigma$, then either $L$ or $N$ is in $\Sigma$.
\end{enumerate}
\end{defi}

Support $\tau$-tilting modules are defined as follows.

\begin{defi}[\cite{AIR}, Definition 0.1]\label{tautiltingmodules}
Let $M$ be a basic $A$-module. Then:
\begin{enumerate}
\item $M$ is said to be \textit{$\tau$-rigid} if $Hom_A(M,\tau_A M)=0$;
\item $M$ is a \textit{$\tau$-tilting} module if $M$ is $\tau$-rigid and $|M|=|A|$;
\item $M$ is a \textit{support $\tau$-tilting} module if there exists an idempotent $e\in A$ such that $M$ is a $\tau$-tilting module over the algebra $A/AeA$.
\end{enumerate} 
\end{defi}

The following proposition records some properties of support $\tau$-tilting modules. 

\begin{prop}[\cite{AIR}, Proposition 2.2]\label{lista}
\begin{enumerate}
\item The $\tau$-tilting modules are precisely the sincere support $\tau$-tilting modules.
\item The tilting modules are precisely the faithful support $\tau$-tilting modules.
\item Any $\tau$-tilting (or $\tau$-rigid) $A$-module T is a tilting (or partial tilting respectively) $C$-module, where $C=A/AnnT$. 
\end{enumerate}
\end{prop}

\begin{prop}\label{propix}
Let $M$ be an $A$-module. Then $M$ is a tilting module if and only if $M$ is a $\tau$-tilting module with pd$M\leq 1$.
\end{prop}

\begin{proof}
Necessity follows from Proposition \ref{lista} (2) and the definition of tilting module. 

For sufficiency, let $M$ be a $\tau$-tilting $A$-module such that $pd_A M\leq 1$. Then $|M|=|A|$ because $M$ is a $\tau$-tilting module. By \cite[Corollary IV.4.7]{ARS}, $Ext^1_A(M,M)\cong DHom_A(M,\tau M)=0$. Therefore $M$ is rigid. Thus $M$ is a tilting module.
\end{proof} 

Let $M$ be a support $\tau$-tilting module. We associate $M$ to the pair $(Fac M,$ $Sub(\tau M))$ of subcategories of mod$A$. Note that if $X\in Fac M$ and $Y\in Sub(\tau M)$ then Hom$_A(X,Y)=0$. Moreover if $M$ is a $\tau$-tilting module, then the pair $(Fac M, Sub(\tau M))$ is a torsion pair, see \cite[Theorem 2.12]{AIR}. 

                                \section{A $\tau$-tilting Theorem}\label{tautiltingtheorem}

Given an algebra $A$ and a support $\tau$-tilting $A$-module $M$, we consider the algebras $B=End_A (M)$ and $C=A/Ann_A M$. We also consider in mod$B$ the torsion pair $(\mathcal{X},\mathcal{Y})$, where $\mathcal{X}=\{X\in \text{ mod}B:X\otimes_BM=0\}$ and $\mathcal{Y}=\{Y\in \text{ mod}B:Tor_1^B(Y,M)=0\}$. The main result of this section is the following.

\begin{theorem}\label{cociente}
Let $M$ be a support $\tau$-tilting $A$-module, $B=End_A M$ be its endomorphism algebra and $C=A/Ann M$. 
\begin{enumerate}
\item The algebra morphism $\varphi: A \longrightarrow End_B(_BM)$ defined by $\varphi(a)(m)=ma$ for every $m\in M$ and $a\in A$ induces an isomorphism $C\cong End_B(_BM)$. 

\item The functor $Ext^1_A(M,-): Sub( \tau_A M)\longrightarrow \mathcal{X}$ is an equivalence of categories with quasi-inverse $Tor_1^B(-,M): \mathcal{X}\longrightarrow Sub( \tau_A M)$ if and only if $\tau_C M=\tau_A M$.
\end{enumerate}
\end{theorem}

\begin{proof}[Proof of 1.]
 Let $M$ be a support $\tau$-tilting $A$-module, $B=End_AM$ and $C=A/AnnM$. Consider the morphism $\varphi: A\rightarrow End_BM$ defined by $a\mapsto (m\mapsto ma)$. For any $a\in A$ we have that $\varphi(a)=0$ if and only if $Ma=0$. Thus $Ker\varphi=Ann M$. Now $\varphi$ induces, by passing to the quotient, a morphism $\tilde{\varphi}:C\rightarrow End_BM$ also defined by $c\mapsto (m\mapsto mc)$.
 $$\xymatrix {
 Ann_A M\ar@ {->} [r] &A \ar@ {->} [d] \ar [rr]^{\varphi} &
 &End_B (_BM) \\ &C\ar@ {-->} [rru]_{\tilde {\varphi}}} $$
 Because of Proposition \ref{lista}, the module $M$ is a $C$-tilting module. Moreover because mod$C$ is a full subcategory of mod$A$, we have that $End_CM\cong End_AM=B$. Therefore it follows from \cite[Lemma VI.3.3]{AsSS} that $\tilde{\varphi}$ is an isomorphism.
\\
\\ 
 \textit{Proof of 2.} Sufficiency. An obvious consequence of the equality of $\tau_C M$ and $\tau_A M$ is that $\tau_A M$ is a $C$-module, so that $Sub(\tau_AM)=Sub(\tau_CM)$. We may then denote $\tau_AM$ and $\tau_CM$ simply by $\tau M$.
 
 We claim that $Ext_A^1(M,-)|_{Sub(\tau M)}=Ext_C^1(M,-)|_{Sub(\tau M)}$. Indeed for $L\in$mod$A$, we have $Ext^1_A(M,L)\cong D\overline{Hom}_A(L,\tau M)$ where $\overline{Hom}_A(L,\tau M)$ is the quotient of $Hom_A(L,\tau M)$ by the subspace $\mathcal{I}_A(L,\tau M)$ of all morphisms from $L$ to $\tau M$ in mod$A$ which factor through an injective $A$-module. Assume $L$ belongs to $Sub(\tau M)$ and $\mathcal{I}_A(L,\tau M)\neq 0$. Then there exist an injective $A$-module $I$ and a nonzero morphism $f:L\rightarrow \tau M$ such that $f$ factors through $I$ as $f=f_2f_1$. As $L$ is a $C$-module, let $\iota: L\rightarrow I'$ be an injective envelope in mod$C$. Then $\iota$ is a monomorphism in mod$A$ and the injectivity of $I$ yields $h:I'\rightarrow I$ in mod$A$ such that $f_1=h\iota$.
 $$\xymatrix{
 L \ar@{->}[dd]_\iota \ar@{->}[rd]_{f_1} \ar[rr]^f
 &
 &\tau M \\
 &I\ar[ru]_{f_2} \\
 I'\ar@{->}_h[ur]}$$
 Therefore $f=f_2h\iota$ and $f_2h\neq 0$ because $f\neq 0$. Moreover, mod$C$ is a full subcategory of mod$A$, and $I'$, $\tau M$ are $C$-modules, so $f_2h$ is a morphism in mod$C$. However, $M$ is a tilting $C$-module, so $pd_C M\leq 1$ and therefore $f_2h=0$, a contradiction. This shows that $\mathcal{I}_A(L,\tau M)=0$. Using again that mod$C$ is a full subcategory in mod$A$, and that $\tau_CM=\tau_AM$ we get functorial isomorphisms
 $$Ext_A^1(M,L)\cong DHom_A(L,\tau M)=DHom_C(L,\tau M)\cong Ext^1_C(M,L)$$
 thus establishing our claim. Because $M$ is a $C$-tilting module, we have that $Ext^1_C(M,-):Sub(\tau M)\rightarrow \mathcal{X}$ is an equivalence with quasi-inverse $Tor_1^B(-,M): \mathcal{X}\longrightarrow Sub(\tau M)$. Hence so is $Ext^1_A(M,-)$.
 
 Necessity. Assume $Ext^1_A(M,-):Sub(\tau M)\rightarrow \mathcal{X}$ is an equivalence of categories. Because $M$ is a tilting $B$-module with $End_B(M)\cong C$, we have that $Sub(\tau_C M)$ is equivalent to $\mathcal{X}$ via the functor $Tor^B_1(-,M)$. Hence $\tau_A M\cong$ $Tor^B_1(Ext^1_A(M,\tau_A M),M)$ is a $C$-module. So it follows from \cite[Lemma VIII.5.2]{AsSS} that $\tau_A M$ is isomorphic to $\tau_C M$.
\end{proof}

As a consequence of the previous result and \cite[Proposition 3.5]{J} we can state a generalization of the tilting theorem of Brenner and Butler.

\begin{cor}\label{tauBB}
Let $M$ be a support $\tau$-tilting $A$-module. Consider the pair of subcategories $(Fac M, Sub(\tau_A M))$ of mod$A$. Let $B=End_A(M)$ be the endomorphism algebra of $M$ and $C= A/Ann M$. Then there exists a torsion pair $(\mathcal{X},\mathcal{Y})$ in mod$B$ such that
\begin{enumerate}
\item The functor $Hom_A(M,-): Fac M\longrightarrow \mathcal{Y}$ is an equivalence of categories with quasi-inverse $-\otimes_B M: \mathcal{Y}\longrightarrow Fac M$.
\end{enumerate}
Moreover $\tau_C M=\tau_A M$ if and only if
\begin{enumerate}
\item[2.] The functor $Ext^1_A(M,-): Sub( \tau_A M)\longrightarrow \mathcal{X}$ is an equivalence of categories with quasi-inverse $Tor_1^B(-,M): \mathcal{X}\longrightarrow Sub( \tau_A M)$.
\end{enumerate}
\end{cor}

\begin{proof}
The first assertion was proved by Jasso in \cite[Proposition 3.5]{J}. The second statement is exactly the second part of Theorem \ref{cociente}.
\end{proof}

It was pointed out by Adachi, Iyama and Reiten in \cite{AIR} that we can develop a dual theory if we consider the $\tau^-$-rigid modules, that is Hom$_A(\tau^{-1}_A M,M)=0$. Following the dual arguments we can prove the following result.

\begin{cor}\label{cotauBB}
Let $M$ be a support $\tau^-$-tilting $A$-module. Consider the pair of subcategories $(Fac \tau_A^{-1}M, Sub M)$ of mod$A$. Let $B=End_A(M)$ be the endomorphism algebra of $M$ and $C= A/Ann M$. Then there exists a torsion pair $(\tilde{\mathcal{X}},\tilde{\mathcal{Y}})$ in mod$B$ such that
\begin{enumerate}
\item The functor $Hom_A(-,M): Sub M\longrightarrow \tilde{\mathcal{X}}$ is an equivalence of categories with quasi-inverse $M\otimes_B -: \tilde{\mathcal{X}}\longrightarrow Sub M$.
\end{enumerate}
Moreover $\tau^{-1}_C M=\tau^{-1}_A M$ if and only if
\begin{enumerate}
\item[2.] The functor $Ext^1_A(-,M): Fac (\tau^{-1}_A M)\longrightarrow \tilde{\mathcal{Y}}$ is an equivalence of categories with quasi-inverse $Tor_1^B(M,-): \tilde{\mathcal{Y}}\longrightarrow Fac (\tau^{-1}_A M)$.
\end{enumerate}
\end{cor}

We now present two examples. The first one shows that $\tau_A M=\tau_C M$ is generally not true for a support $\tau$-tilting module $M$. The second example shows a support $\tau$-tilting module $M$ such that $\tau_AM=\tau_CM$. 

\begin{ex}
Let $A$ be the path algebra of the quiver
$$\xymatrix{1 \ar@/^/[r]^\alpha & 2 \ar@/^/[l]^{\alpha '}\ar@/^/[r]^\beta & 3 \ar@/^/[l]^{\beta '}}$$
bound by the ideal $I=\langle\alpha'\alpha-\beta\beta', \alpha\alpha', \beta'\beta\rangle$. Take the $\tau$-tilting module $M=\rep{1\\2\\3}\oplus \rep{1\\2} \oplus \rep{1}$. Then $AnnM=\langle\alpha', \beta'\rangle$ and the quotient algebra $C=A/AnnM$ is the path algebra of the following quiver of type $\mathbb{A}_3$.
$$\xymatrix{1\ar[r]^{\alpha}&2\ar[r]^\beta&3}$$
We know that $M$ is a tilting module in mod$C$. Hence there exists a torsion pair $(\mathcal{X},\mathcal{Y})$ in mod$B$, where $B=End_A(M)=End_C(M)$. This torsion pair is such that the functor Hom$_C(M,-)$ is an equivalence of categories between Fac$M$ and $\mathcal{Y}$; while Ext$^1_C(M,-)$ is an equivalence of categories between $Sub(\tau_C M)$ and $\mathcal{X}$. Now $\tau_CM=\rep{2\\3}\oplus \rep{2}$ and $\tau_AM=\rep{2\\3}\oplus \rep{3\\2}$. This implies that $Sub(\tau_C M)$ is strictly contained in $Sub(\tau_A M)$. So the torsion class $\mathcal{X}$ in mod$B$ is not equivalent to $Sub(\tau_A M)$.
\end{ex}

\begin{ex}
Consider the algebra $A=kQ/I$, where $Q$ is
$$
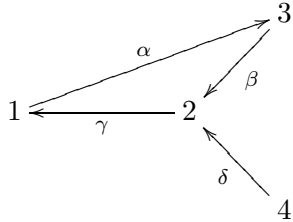
\xymatrix{
 & & & 3\ar[dl]^{\beta} \\
1\ar[urrr]^{\alpha} & & 2\ar[ll]^{\gamma} & \\
 & & &4\ar[ul]^{\delta} }
$$
and $I=\langle \alpha\beta, \gamma\alpha\rangle$. Take $M=\rep{4\\2\\1}\oplus\rep{4\\2}\oplus\rep{43\\2\\1}\oplus\rep{4}$, with Auslander-Reiten translate $\tau_{A}M=\rep{3\\2\\1}\oplus\rep{2}\oplus\rep{3\\2}$. Then Hom$_{A}(M, \tau_{A}M)=0$ and $|M|=|A|=4$. Therefore $M$ is a $\tau$-tilting $A$-module. Moreover $M$ is not a tilting module because its annihilator $Ann_{A} M=\langle \alpha\rangle$ is not zero. Consider the torsion pair $(Fac M, Sub(\tau_A M))$ associated to $M$, then $Fac M= add\{\rep{4\\2\\1}, \rep{4\\2}, \rep{43\\2\\1}, \rep{4}, \rep{43\\2}, \rep{3}\}$ and $Sub (\tau_{A}M)=add\{\rep{3\\2\\1}, \rep{2}, \rep{3\\2}, \rep{2\\1}, \rep{1}\}$.

Consider the algebra $C=A/Ann_{A}M$. Then $C$ is the hereditary algebra given by the following quiver.
$$
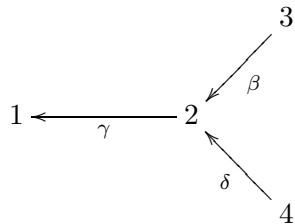
\xymatrix{
 & & & 3\ar[dl]^{\beta} \\
1 & & 2\ar[ll]^{\gamma} & \\
 & & &4\ar[ul]^{\delta} }
$$
We have that $\tau_{C}M=\rep{3\\2\\1}\oplus\rep{2}\oplus\rep{3\\2}\cong\tau_{A}M$.
On the other hand, $B=End_{A}(M)=kQ'/I'$, where $Q'$ is the quiver
$$
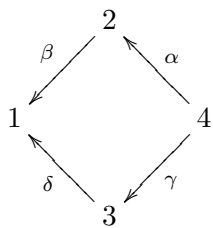
\xymatrix{
 &2\ar[ld]_{\beta} & \\
 1& &4\ar[lu]_{\alpha}\ar[ld]^{\gamma} \\
 &3\ar[lu]^{\delta} &
}
$$
and $I'$ is the ideal generated by the commutativity relation $\alpha\beta - \gamma\delta$. Therefore, Corollary \ref{tauBB} yields a torsion pair $(\mathcal{X}, \mathcal{Y})$ in mod$B$ such that
$$Hom_A(M,-): Fac M\longrightarrow \mathcal{Y}$$
$$Ext^1_A(M,-): Sub( \tau_A M)\longrightarrow \mathcal{X}$$
are equivalences of categories. In this particular case we can calculate explicitly the pair $(\mathcal{X},\mathcal{Y})$ to get $\mathcal{X}=add\{\rep{3}, \rep{4\\32}, \rep{4\\3}, \rep{4\\2}, \rep{4}\}$ and $\mathcal{Y}=add\{ \rep{4\\32\\1}, \rep{2}, \rep{32\\1}, \rep{2\\1}, \rep{3\\1}, \rep{1}\}$.
\end{ex}

Two other consequences of Theorem \ref{cociente} are the following. 

\begin{cor}
Let $M$ be a support $\tau$-tilting $A$-module and $B=End_A(M)$ be its endomorphism algebra. Then $B$ is connected if and only if $C=A/Ann M$ is connected.
\end{cor}

\begin{proof}
Let $M$ be a support $\tau$-tilting $A$-module. Then $M$ is a tilting $C$-module. Therefore the center $Z(C)$ of the algebra $C$ is isomorphic to the center $Z(B)$ of the algebra $B$ (see \cite[Lemma VI.3.4]{AsSS}). Hence $B$ is connected if and only if $C$ is connected.
\end{proof}

\begin{cor}
Let $M$ be a support $\tau$-tilting $A$-module and $B=End_A(M)$ be its endomorphism algebra. Then the algebra $End_B(M)$ is isomorphic to the algebra $A$ if and only if $M$ is a tilting module.
\end{cor}

\begin{proof}
Consider the morphism $\varphi$ from Theorem \ref{cociente} and suppose that $\varphi$ is an isomorphism of algebras. Then it follows from Theorem \ref{cociente} that $Ann_A M=0$. Consequently we can apply Proposition \ref{lista} to conclude that $M$ is an tilting $A$-module. The other implication is shown in \cite[Lemma VI.3.3]{AsSS}.
\end{proof}

\section{$\tau$-slices}

There is a particular class of $\tau$-tilting modules which behaves much like the complete slices in tilted algebras. We call them \textit{$\tau$-slices}. The rest of the paper is devoted to their study.

\subsection{Definition}

The $\tau$-slices are support $\tau$-tilting $A$-modules whose indecomposable direct summands induce a presection in the Auslander-Reiten quiver of $A$.

\begin{defi}[\cite{ABS2}, Definition 3]\label{presection}
Let $A$ be a finite dimensional $k$-algebra, $\Gamma_A$ be its Auslander-Reiten quiver and $\Gamma\subset\Gamma_A$ be a connected component. A full connected subquiver $\Sigma$ of $\Gamma$ is called a \textit{presection} if:
\begin{enumerate}
\item for any given arrow $f:X\rightarrow Y$ in $\Gamma$ such that $X\in\Sigma$, either $Y$ or $\tau Y$ belongs to $\Sigma$.
\item for any given arrow $f:X\rightarrow Y$ in $\Gamma$ such that $Y\in\Sigma$, either $X$ or $\tau^{-1} X$ belongs to $\Sigma$.
\end{enumerate}
\end{defi}

The following are consequences of the definition of presection.

\begin{lem}\label{presecpro}
Let $\Sigma$ be a presection. Then no immediate sucessor of $\Sigma$ is projective. Dually, no immediate predecessor of $\Sigma$ is injective.
\end{lem}

\begin{proof}
Let $P$ be an indecomposable projective module such that there exists an indecomposable summand $M$ of rad$P$ which belongs to $\Sigma$. Then there exists an arrow from $M$ to $P$ in $\Gamma$. Since $\Sigma$ is a presection, either $P$ or $\tau P$ belongs to $\Sigma$. But $\tau P = 0$, then $P$ belongs to $\Sigma$. 
The other case is dual. 
\end{proof}

\begin{lem}\label{presecpro2}
Let $\Sigma$ be a presection. Then $\Sigma$ is acyclic.
\end{lem}

\begin{proof}
By definition of presection, every path contained in $\Sigma$ is a sectional path, see \cite[Section IX.2]{AsSS}. Therefore \cite[Corollary IX.2.3]{AsSS} implies that $\Sigma$ is acyclic.
\end{proof}

Now we are able to give the formal definition of the $\tau$-slices. 

\begin{defi}\label{tauslice}
Let $\Sigma$ be a presection in the Auslander-Reiten quiver $\Gamma_A$ of $A$. Then $\Sigma$ is a \textit{$\tau$-slice} if $M=\bigoplus\limits_{U\in\Sigma}U$ is a support $\tau$-tilting module. It is a \textit{complete $\tau$-slice} if $M$ is a $\tau$-tilting module. 
\end{defi}

\begin{rmk}
If $\Sigma$ is a complete $\tau$-slice, then it follows from Definition \ref{tautiltingmodules} that $|\Sigma|=|A|$. Moreover $\Sigma$ is sincere by Proposition \ref{lista}.
\end{rmk}

\begin{rmk}
From now on, by abuse of notation, we denote by $\Sigma$ the presection in $\Gamma_A$ and the corresponding $A$-module. The context will indicate if we are speaking about the module in the category mod$A$ or the subquiver of $\Gamma_A$.
\end{rmk}

One key property of $\tau$-slices is the following.

\begin{prop}\label{endoslice}
Let $\Sigma$ be a $\tau$-slice in mod$A$. Then the algebra $B=$End$_A(\Sigma)$ is a hereditary algebra. 
\end{prop}

\begin{proof}
Let $P$ be an indecomposable projective $B$-module and let $u:Q\rightarrow P$ be a nonzero monomorphism in mod$B$. Consider the torsion pair $(\mathcal{X},\mathcal{Y})$ associated to $_B\Sigma$ viewed as a left $B$-module. Because $\Sigma$ is a support $\tau$-tilting module, \cite[Proposition 3.5]{J} implies that Hom$_A(\Sigma, -)$ induces an equivalence of categories between $Fac_A \Sigma$ and $\mathcal{Y}$. Since $P$ belongs to $\mathcal{Y}$ so does $Q$ because $\mathcal{Y}$ is closed under submodules. In addition, \cite[Proposition 3.5]{J} implies the existence of $M\in\Sigma$, $N\in Fac\Sigma$ and a nonzero morphism $v: N\rightarrow M$ such that Hom$_A(\Sigma,v)=u$. Note that $\Sigma$ induces a finite presection $\Sigma$ in $\Gamma_A$, which is acyclic by Lemma \ref{presecpro2}. If $N$ is not a direct summand of $\Sigma$, then \cite[Lemma VIII.5.4]{AsSS} yields a factorization $v=gf$ where $f:N\rightarrow \tau_A \Sigma$ and $g:\tau_A \Sigma \rightarrow \Sigma$. But Hom$_A(\Sigma, \tau_A\Sigma)=0$. Therefore $f=0$. This contradiction emerged from the assumption that $N$ does not belong to $\Sigma$. Therefore $N\in\Sigma$, which implies that $Q$ is a projective $B$-module. 
\end{proof}

The following result first appeared in \cite{Liu} with an equivalent formulation. We give here an alternative proof, using the results developed in Section \ref{tautiltingtheorem}.

\begin{cor}\label{tiltedcociente}
Let $\Sigma$ be a $\tau$-slice in mod$A$. Then $C=A/Ann \Sigma$ is a tilted algebra. Moreover $\Sigma$ is a complete slice in mod$C$.
\end{cor}

\begin{proof}
Because $\Sigma$ is a support $\tau$-tilting module, it follows from Theorem \ref{cociente} that $C\cong End_B(_B\Sigma)$, where $B=End_A(\Sigma)$. Moreover $B$ is hereditary by Proposition \ref{endoslice}. Therefore $C$ is the endomorphism algebra of a tilting module over a hereditary algebra, and thus is tilted.
Finally, $\Sigma_C\cong Hom_A(\Sigma,\Sigma)\otimes_B\Sigma\cong B\otimes_B\Sigma$ which is a complete slice in mod$C$ because $\Sigma$ is a tilting $B$-module and $B$ is hereditary, see \cite[Corollary 4.2.3]{Rin}.
\end{proof}

The next lemma will be useful subsequently.

\begin{lem}\label{rigidezalcanza}
Let $\Sigma$ be a presection in $\Gamma_A$ such that $\Sigma$, as an $A$-module, is $\tau$-rigid. Then $\Sigma$ is a $\tau$-slice in $\Gamma_A$.
\end{lem}

\begin{proof}
Applying Corollary \ref{tiltedcociente} we know that $\Sigma$ is a complete slice in the tilted algebra $C=A/Ann\Sigma$. We have that $|\Sigma|=|C|$ because $\Sigma$ is a tilting $C$-module. Take the set of primitive idempotents $S := \{e_i : e_i \in Ann \Sigma\}$, the idempotent $e = \sum\limits_{e_i \in S} e_i$ and the algebra $A'=A/ AeA$. Then the annihilator of $\Sigma$ over $A'$ is contained in rad$A'$ because if $e\in A'$ is an idempotent such that $\Sigma e=0$, then $e=0$. Thus $|C|=|A'|$. Therefore $\Sigma$ is a $\tau$-rigid $A'$-module such that $|\Sigma|=|A'|$, which means that $\Sigma$ is a support $\tau$-tilting $A$-module. As $\Sigma$ is also a presection by hypothesis, $\Sigma$ is a $\tau$-slice.
\end{proof}

It is known that complete slices are convex inside the module category. Following Liu \cite{Liu}, we consider a weaker notion of convexity.

\begin{cor}\label{wconvex}
Let $\Sigma$ be a $\tau$-slice in mod$A$. Then $\Sigma$ is weakly convex, that is, for every path 
$$X\overset{f_0}\longrightarrow X_1\overset{f_1}\longrightarrow\dots\overset{f_{n-1}}\longrightarrow X_n\overset{f_n}\longrightarrow Y$$
in $\Gamma_A$ such that the morphism $f=f_n...f_0$ is different from zero and $X,Y\in\Sigma$, we have $X_i\in\Sigma$ for all $i\in\{1,...,n\}$.
\end{cor}

\begin{proof}
 The $\tau$-slice $\Sigma$ is a $\tau$-rigid presection. Therefore $\Sigma$ is weakly convex by \cite[Proposition 2.5]{Liu}.
\end{proof}

        \subsection{Quotients of Algebras with $\tau$-Slices}

Let $A$ be an algebra having a $\tau$-slice $\Sigma$. In this subsection, we study the quotient algebras $A/I$, where $I\subseteq Ann\Sigma$. As a consequence of the main results of this subsection we get that $\tau_A \Sigma=\tau_C \Sigma$, where $C=A/Ann\Sigma$, and thus Corollary \ref{tauBB} applies.

The next two propositions are needed in the proof of the main theorem of this subsection.

\begin{lem}\label{prop1}
Let $I$ be an ideal of $A$ such that $I\subseteq Ann\Sigma$. If $f:X\rightarrow Y$ is an irreducible morphism, where $X\in\Sigma$, then $Y$ is an $A/I$-module. Dually, if $Y\in\Sigma$, then $X$ is an $A/I$-module.
\end{lem}

\begin{proof}
We prove the case where $X$ belongs to $\Sigma$. The other is dual. 

First notice that $\Sigma$ is an $A/I$-module because $I\subset Ann\Sigma$. Then every indecomposable summand of $\Sigma$ is also an $A/I$-module. Suppose that $X$ belongs to $\Sigma$ and that $Y$ does not. Then Lemma \ref{presecpro} implies that $Y$ is not a projective $A$-module. Hence $\tau_A Y$ belong to $\Sigma$ because $\Sigma$ is a presection in $\Gamma_A$. Consider the almost split sequence ending with $Y$.
$$0\rightarrow\tau_A Y\overset{g}\longrightarrow E\overset{g'}\longrightarrow Y\rightarrow 0$$
If $E$ is itself an $A/I$-module, then so does $Y$ because it is the cokernel of a morphism in mod$A/I$. Suppose that $E$ is not an $A/I$-module. Then there exists an indecomposable direct summand $Y_1$ of $E$ such that $Y_1$ does not belong to mod$A/I$. So, there exists an arrow $f_1$ in $\Gamma_A$ such that $f_1:\tau_A Y\rightarrow Y_1$. We know that $\tau_A Y$ belong to $\Sigma$. Fix an irreducible morphism $h_1:Y_1\rightarrow Y$ from $Y_1$ to $Y$.

Consider the almost split sequence ending with $Y_1$.
$$0\rightarrow\tau_A Y_1\overset{g_1}\longrightarrow E_1\overset{g'_1}\longrightarrow Y_1\rightarrow 0$$
Then $\tau_A Y_1$ belongs to $\Sigma$ because $\tau_A Y$ is a direct summand of $E_1$. Also there exists a direct summand $Y_2$ of $E_1$ which does not belong to mod$A/I$. Fix an irreducible morphism $h_2:Y_2\rightarrow Y_1$ from $Y_2$ to $Y_1$.

Inductively, for every natural number $n$ we have a module $Y_n$ which is not an $A/I$-module and such that $\tau_A Y_n$ is a direct summand of $\Sigma$ and an irreducible morphism $h_n:Y_n\rightarrow Y_{n-1}$. The composition $h_nh_{n-1}\dots h_2h_1$ is a sectional path for every $n\in\mathbb{N}$. Otherwise there exists $i\leq n$ such that $Y_i=\tau_A Y_{i-2}$ which belongs to $\Sigma$, a contradiction.

Thus, we have an infinite sectional path in $\Gamma_A$. Then $Y_i\neq Y_j$ for $i\neq j$ by \cite[Corollary IX.2.3]{AsSS}. Moreover $\tau_A Y_i\in\Sigma$ for every $i$, which contradicts the finiteness of $\Sigma$. Therefore $Y$ is an $A/I$-module.
\end{proof}

\begin{lem}\label{prop2}
Let $I$ be an ideal of $A$ such that $I\subseteq Ann\Sigma$ and $f:X\rightarrow Y$ be a morphism with $X,Y$ indecomposable. Assume $X$ or $Y$ belongs to $\Sigma$. Then $f$ is irreducible in mod$A$ if and only if $f$ is irreducible in mod$A/I$.
\end{lem}

\begin{proof}
Suppose $X$ is in $\Sigma$, the case where $Y$ is in $\Sigma$ is dual. 

The necessity follows from Lemma \ref{prop1}, so we prove the sufficiency. Let $f$ be irreducible in mod$A/I$. In particular both $X$ and $Y$ are $A/I$-modules. So $f$ is a morphism in mod$A$ but not necessarily irreducible. Let $X\overset{\left(\rep{f_1\\ \vdots \\ f_t}\right)}\longrightarrow E_1\oplus\dots\oplus E_t$ be left almost split in mod$A$, with $E_i$ indecomposable for all $i$. Then $f$ is not a section in mod$A$ because it is not a section in mod$A/I$. Then there exists $[f'_1, \dots, f'_t]:E_1\oplus\dots\oplus E_t\longrightarrow Y$ such that $f=\sum_{i=1}^t f'_if_i$. Lemma \ref{prop1} implies that $E_i$ is an $A/I$-module for every $i$. Therefore the factorization given in the sum is a factorization of $f$ in mod$A/I$. By hypothesis $f$ is irreducible in mod$A/I$. Therefore $[f_1, \dots, f_t]^t$ is a section or $[f'_1, \dots, f'_t]$ is a retraction. Because the $f_i$ are irreducible, one of the $f'_i$ is a retraction and then an isomorphism. Hence $f$ is irreducible in mod$A$.
\end{proof}

\begin{cor}\label{tauiguales}
Let $I$ be an ideal of $A$ such that $I\subseteq Ann\Sigma$. Then $\tau_{A/I}X=\tau_AX$ and $\tau^{-1}_{A/I}X=\tau^{-1}_AX$ for every indecomposable direct summand $X$ of $\Sigma$.
\end{cor}

\begin{proof}
 Let $X$ be a indecomposable summand of $\Sigma$ and consider the almost split short exact sequence in mod$A$ ending with $X$.
$$0\rightarrow\tau_A X\overset{u}\longrightarrow M\overset{v}\longrightarrow X\rightarrow 0$$

By Lemmata \ref{prop1} and \ref{prop2} the almost split sequence ending with $X$ in mod$A$ coincides with the almost split sequence ending with $X$ in mod$A/I$. Hence $\tau_A X=\tau_{A/I} X$. Dually we can prove that $\tau^{-1}_AX=\tau^{-1}_{A/I}X$.
\end{proof}

Now we can prove the main theorem of this section.

\begin{theorem}\label{teorema1}
Let $\Sigma$ be a $\tau$-slice in mod$A$ and let $I$ be an ideal of $A$ such that $I\subseteq Ann \Sigma$. Then $\Sigma$ is also a $\tau$-slice in mod$A/I$.
\end{theorem}

\begin{proof}
First we prove that $\Sigma$ is a presection in mod$A/I$. Let $f:X\rightarrow Y$ be an arrow in $\Gamma_{A/I}$ such that $X\in\Sigma$. If $Y\in\Sigma$, there is nothing to prove. Suppose the contrary. Now, $f$ is an irreducible morphism in mod$A$ which is a direct summand of the right almost split morphism $f':E\rightarrow Y$. So, we can write the almost split sequence ending with $Y$. 
$$0\rightarrow\tau_A Y\overset{g'}\longrightarrow E\overset{f'}\longrightarrow Y\rightarrow 0$$
with $X$ an indecomposable direct summand of $E$. 

Then, because $\Sigma$ is a presection in mod$A$, it follows that $\tau_AY\in\Sigma$. Then Corollary \ref{tauiguales} implies that $Y=\tau^{-1}_A \tau_A Y=\tau^{-1}_{A/I}\tau_A Y$. Hence $\tau_A Y=\tau_{A/I} Y$. Therefore $\tau_{A/I} Y\in\Sigma$. 

Dualizing arguments, we prove that if $f:X\rightarrow Y$ is an arrow in $\Gamma_A$ such that $Y\in\Sigma$ then either $X$ or $\tau^{-1}_{A/I}X$ belongs to $\Sigma$. Therefore $\Sigma$ is a presection in mod$A/I$. 

Moreover $Hom_{A/I}(\Sigma,\tau_{A/I}\Sigma)=Hom_A(\Sigma,\tau_{A/I}\Sigma)=Hom_A(\Sigma,\tau_A\Sigma)=0$. Therefore $\Sigma$ is a $\tau$-rigid presection in mod$A/I$. Lemma \ref{rigidezalcanza} implies that $\Sigma$ is a $\tau$-slice in mod$A/I$, finishing the proof.
\end{proof}

\begin{rmk}\label{rmk1}
Let $\Sigma$ be a $\tau$-slice in mod$A$. Then $\Sigma$ is a complete slice in mod$C$ by Corollary \ref{tiltedcociente}. Therefore \cite[Lemma VIII.5.5]{AsSS} implies that Hom$_C(\tau^{-1}_C \Sigma,\Sigma)=0$. Hence Hom$_A(\tau^{-1}_A\Sigma,\Sigma)=0$ by Corollary \ref{tauiguales}. Thus every $\tau$-slice $\Sigma$ is at the same time a support $\tau$-tilting module and a support $\tau^{-1}$-cotilting module in mod $A$. This behavior of $\tau$-slices generalizes the property of the complete slices which are tilting modules and cotilting modules at the same time.
\end{rmk}

\begin{cor}\label{slicesBB}
Let $\Sigma$ be a $\tau$-slice in mod$A$ and consider the torsion pair $(Fac \Sigma, Sub(\tau \Sigma))$. Let $B=End_A(\Sigma)$ be the endomorphism algebra of $\Sigma$. Then there exists a torsion pair $(\mathcal{X},\mathcal{Y})$ in mod$B$ such that
\begin{enumerate}
\item The functor $Hom_A(\Sigma,-): Fac \Sigma\longrightarrow \mathcal{Y}$ is an equivalence of categories with quasi-inverse $-\otimes_B\Sigma: \mathcal{Y}\longrightarrow Fac \Sigma$.
\item The functor $Ext^1_A(\Sigma,-):  Sub(\tau \Sigma)\longrightarrow \mathcal{X}$ is an equivalence of categories with quasi-inverse $Tor_1^B(-,\Sigma): \mathcal{X}\longrightarrow  Sub(\tau \Sigma)$.
\end{enumerate}
\end{cor}

\begin{proof}
Theorem \ref{teorema1} tells us that $\tau_A \Sigma \cong \tau_{A/I} \Sigma$ for every ideal $I\subseteq Ann_A\Sigma$. In particular $\tau_A \Sigma \cong \tau_{C} \Sigma$, where $C=A/Ann_A \Sigma$. Then the statement follows directly from Corollary \ref{tauBB}.
\end{proof}

\subsection{Tilted and Cluster Tilted Algebras}

In this subsection, we study complete slices in tilted algebras and local slices in cluster tilted algebras as $\tau$-slices.

\begin{lem}\label{completeslices}
Let $A$ be a tilted algebra and $\Sigma$ be a complete slice in mod$A$. Then $\Sigma$ is a complete $\tau$-slice in mod$A$.
\end{lem}

\begin{proof}
First, it was shown in \cite{HR} that every complete slice $\Sigma$ is a tilting module. Then $\Sigma$ is a faithful $\tau$-tilting module. Moreover we know by \cite{LiuCrit} and \cite{Skocrit} that $\Sigma$ induces a section in $\Gamma_A$. Then $\Sigma$ is a presection by \cite[Lemma VIII.1.4]{AsSS}. Thus Lemma \ref{rigidezalcanza} implies that $\Sigma$ is a complete $\tau$-slice.
\end{proof}

The following characterization of tilted algebras is a consequence of our results so far. It was first shown by Liu \cite{Liu} and independently by Skowronski and Yamagata (private communication).

\begin{cor}[\cite{Liu}, Theorem 2.6]\label{taurodajasfieles}
An algebra $A$ is tilted if and only if there exists a faithful $\tau$-slice $\Sigma$ in mod$A$.
\end{cor}

\begin{proof}
Let $A$ be a tilted algebra. Then there exists a complete slice $\Sigma$ in mod$A$. Lemma \ref{completeslices} implies that $\Sigma$ is a complete $\tau$-slice. Moreover $\Sigma$ is faithful because it is a tilting $A$-module. 

Let $A$ be an algebra with a faithful $\tau$-slice $\Sigma$. Because $\Sigma$ is faithful, then $Ann_A\Sigma=0$. Hence Corollary \ref{tiltedcociente} implies that $A$ is a tilted algebra.
\end{proof}

Another type of slice present in the literature is the local slices defined by Assem, Br\"ustle and Schiffler in \cite{ABS2}.

\begin{defi}[\cite{ABS2}, Definition 11]
Let $A$ be an algebra. A subset $\Sigma$ of the Auslander-Reiten quiver $\Gamma_A$ of $A$ is a \textit{local slice} if the following conditions hold.
\begin{enumerate}
\item $\Sigma$ is a presection.
\item $\Sigma$ is \textit{sectionally convex}, that is for every path 
$$X=X_0\overset{f_0}\longrightarrow X_1\overset{f_1}\longrightarrow\dots\overset{f_{n-1}}\longrightarrow X_n=Y$$
in $\Gamma_A$ such that $X_i\neq \tau X_{i+2}$ for $0\leq i\leq n-2$ and $X,Y\in\Sigma$, we have $X_i\in\Sigma$ for all $i\in\{1,...,n-2\}$.
\item $|\Sigma|=|A|$
\end{enumerate}
\end{defi}

\begin{lem}\label{lemix}
Let $\Sigma$ be a $\tau$-slice. Then $\Sigma$ is sectionally convex. 
\end{lem}

\begin{proof}
Consider the sectional path $X=X_0\overset{f_0}\longrightarrow X_1\overset{f_1}\longrightarrow\dots\overset{f_{n-1}}\longrightarrow X_n=Y$, such that $X,Y\in\Sigma$. Then the map $f=f_{n-1}...f_1\neq 0$ (see \cite[Corollary IX.2.2]{AsSS}). Therefore $X_i\in\Sigma$ because $\Sigma$ is weakly convex by Corollary \ref{wconvex}. Thus $\Sigma$ is sectionally convex. 
\end{proof}

The following proposition shows that the notions of $\tau$-slices and local slices coincide in cluster tilted algebras.

\begin{prop}\label{clustertilted}
Let $A$ be an algebra and $\Sigma$ a complete $\tau$-slice. Then $\Sigma$ is a local slice. Moreover, if $A$ is a cluster-tilted algebra then $\Sigma$ is a local slice if and only if $\Sigma$ is a complete $\tau$-slice.
\end{prop}

\begin{proof}
Let $\Sigma$ be a complete $\tau$-slice in mod$A$. By definition $\Sigma$ is a presection in $\Gamma_A$. Moreover $|\Sigma|=|A|$ and $\Sigma$ is sectionally convex by Lemma \ref{lemix}. Therefore $\Sigma$ is a local slice.

Let $A$ be a cluster-tilted algebra and let $\Sigma$ be a local slice in $\Gamma_A$. Therefore $\Sigma$ is a presection and $|A|=|\Sigma|$. Let $\mathcal{C}$ be the cluster category and $T$ be a cluster-tilting object in $\mathcal{C}$ such that $A\cong End_{\mathcal{C}}(T)$.  It was shown in \cite{ABS2} that for every local slice $\Sigma$ in $\Gamma_A$ there exists a cluster-tilting object $\tilde{\Sigma}$ such that $\Sigma=Hom_{\mathcal{C}}(T,\tilde{\Sigma})$. Without loss of generality we can suppose that $\tilde{\Sigma}$ and $\tau_{\mathcal{C}}T$ do not have direct summands in common. Now, \cite[Theorem A]{BMR} implies that $Hom_A(\Sigma,\tau_A \Sigma)\cong Hom_A(Hom_{\mathcal{C}}(T,\tilde{\Sigma}), Hom_{\mathcal{C}}(T,\tau_{\mathcal{C}}\tilde{\Sigma}))\cong Hom_{\mathcal{C}}(\tilde{\Sigma},\tau_{\mathcal{C}}\tilde{\Sigma})/\tilde{\mathcal{I}}$, where $\tilde{\mathcal{I}}$ is the set of morphisms from $\tilde{\Sigma}$ to $\tau_{\mathcal{C}}\tilde{\Sigma}$ that factor through an object in the category $add(\tau_{\mathcal{C}}T)$. But $Hom_{\mathcal{C}}(\tilde{\Sigma}, \tau_{\mathcal{C}}\tilde{\Sigma})=Hom_{\mathcal{C}}(\tilde{\Sigma}, \tilde{\Sigma}[1])=0$ because $\tilde{\Sigma}$ is cluster tilting in $\mathcal{C}$. Therefore $Hom_A(\Sigma,\tau_A\Sigma)=0$. It follows from Lemma \ref{rigidezalcanza} that $\Sigma$ is a complete $\tau$-slice. 
\end{proof}

\begin{question}
Let $A$ be an algebra and $\Sigma$ be a local slice in mod$A$. Is $\Sigma$ a complete $\tau$-slice in mod$A$?
\end{question}

Recall from Proposition \ref{endoslice} that if $\Sigma$ is a $\tau$-slice in mod$A$, then End$_A(\Sigma)$ is hereditary. Given a $\tau$-slice $\Sigma$ we say that $\Sigma$ is of \textit{type $\Delta$} if End$_A(\Sigma)=k\Delta$. The following example show us an interesting fact: if $\Sigma_1$ and $\Sigma_2$ are complete $\tau$-slices in mod$A$, then the type of $\Sigma_1$ might be different from the type of $\Sigma_2$.

\begin{ex}
Take the path algebra $A$ of the following quiver 

$$
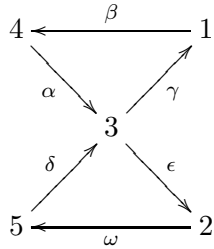
\xymatrix{
4 \ar@{->}[rd]_{\alpha} 
&
&1 \ar[ll]_\beta\\
&3 \ar[ru]_{\gamma} \ar[rd]^\epsilon \\
5\ar@{->}^\delta[ur]
&
&2\ar[ll]^\omega
}$$
bound by the ideal generated by all paths of length two except $\alpha\epsilon$. Its Auslander-Reiten quiver is drawn in Figure 1.

\begin{figure}\label{figura1}
$$
\xymatrix{
 \rep{2\\5}\ar@{->}[rd]& & & & & & & &  & \\ 
  & \rep{2}\ar@{->}[rd]& & \rep{3\\1}\ar@{->}[rd] & & & & & & \\ 
  & & \rep{3\\21}\ar@{->}[rd]\ar@{->}[ru] & & \rep{3}\ar@{->}[rd]\ar@{->}[r]& \rep{5\\3}\ar@{->}[r] & \rep{54\\3}\ar@{->}[r]\ar@{->}[rd]&\rep{4}\ar@{->}[r] &\rep{1\\4}\ar@{->}[r] &\rep{1} \\  
  & \rep{1}\ar@{->}[ru] & & \rep{3\\2}\ar@{->}[rd]\ar@{->}[ru] & & \rep{4\\3}\ar@{->}[ru] & & \rep{5}\ar@{->}[rd] & & \rep{2} \\   
 \rep{1\\4}\ar@{->}[ru] & & & & \rep{4\\3\\2}\ar@{->}[ru] & & & & \rep{2\\5}\ar@{->}[ru]& 
}$$
    \caption{}
\end{figure}

One can see that $\Sigma=\rep{4\\3\\2}\oplus \rep{4\\3}\oplus \rep{54\\3}\oplus \rep{4}\oplus \rep{1\\4}$ is a complete $\tau$-slice of type $\mathbb{A}_5$ while $\tilde{\Sigma}= \rep{4\\3\\2}\oplus \rep{4\\3}\oplus \rep{3}\oplus \rep{5\\3}\oplus \rep{3\\1}$ is a complete $\tau$-slice of type $\mathbb{D}_4$.
\end{ex}

\section{$\tau$-slices and Tilted Algebras}

In this section we study complete $\tau$-slices in tilted algebras. We now prove the key lemma of this section.

\begin{lem}\label{caminito}
Let $\Gamma$ be a connected component of the Auslander-Reiten quiver $\Gamma_A$ of an algebra $A$. Assume that $\Gamma$ is convex, generalized standard and has a complete $\tau$-slice $\Sigma$.

If there is no path in $\Gamma$ from $\Sigma$ to a projective and, dually, no path in $\Gamma$ from an injective to $\Sigma$ then $A$ is tilted having $\Sigma$ as complete slice. 
\end{lem}

\begin{proof}
Let $C=A/Ann\Sigma$. By Corollary \ref{tiltedcociente}, $C$ is tilted having $\Sigma$ as complete slice. We prove that the component $\Gamma$ consists only of $C$-modules. To do so, we show by induction that every predecessor $X$ of $\Sigma$ is the kernel of a morphism in mod$C$, and thus a $C$-module. The proof that every successor of $\Sigma$ is the cokernel of a morphism in mod$C$ is dual and is left to the reader.

Because $\Sigma$ is a finite presection, $\Sigma$ has at least one sink $X$ by Lemma \ref{presecpro2}. If $X$ is projective then $\tau X=0$ and trivially it is a $C$-module. Suppose that $X$ is not projective. Since $X$ is a sink in $\Sigma$ then, for each arrow $f:Y\rightarrow X$ in $\Gamma$, the module $Y$ belongs to $\Sigma$. Take the following almost split sequence.
$$0\rightarrow\tau_A X\rightarrow Y\rightarrow X\rightarrow 0$$
Then every indecomposable summand of $Y$ is a $C$-module. Consequently, so is $Y$. Therefore $\tau_A X$ is a $C$-module because it is the kernel of a morphism in mod$C$.

Let $X$ be a predecessor of $\Sigma$, such that, inductively, $\tau^{-1}_A X$ and every immediate successor of $X$ is a $C$-module. By hypothesis $X$ is not an injective $A$-module, then $\tau^{-1}_AX\neq 0$. Consider the almost split sequence starting with $X$.
$$0\longrightarrow X\overset{g}\longrightarrow Y\overset{f}\longrightarrow \tau^{-1}_AX\longrightarrow 0$$

Then $Y$ and $\tau^{-1}_AX$ are $C$-modules. Therefore $f$ is a morphism in mod$C$. Thus $X$ is a $C$-module. This completes the proof that the component $\Gamma$ consists of $C$-modules.

We know that $\Sigma$ is a complete slice of mod$C$. Then $\Sigma$ is convex in mod$C$. In particular $\Sigma$ is convex in $\Gamma$. Moreover $\Gamma$ is convex in mod$A$. Therefore $\Sigma$ is convex in mod$A$.

Because $\Sigma$ is a complete $\tau$-slice, $\Sigma$ is sincere as $A$-module. Take an almost split sequence 
$$0\longrightarrow X\longrightarrow Z\longrightarrow Y\longrightarrow 0$$
with $X$ and $Y$ indecomposable. Suppose that an indecomposable summand $Z_1$ of $Z$ belongs to $\Sigma$. Then there is an arrow $f:X\rightarrow Z_1$ from $X$ to $Z_1$. Therefore either $X$ or $Y$ belongs to $\Sigma$ because $\Sigma$ is a presection.

Then $\Sigma$ is a complete slice in mod$A$ and the statement follows.
\end{proof}

The following theorem is a consequence of Lemma \ref{caminito}.

\begin{theorem}\label{taurigidsection}
Let $\Gamma$ be a connected component of the Auslander-Reiten quiver $\Gamma_A$ of an algebra $A$. Assume that $\Gamma$ is convex and generalized standard. Then $A$ is tilted if and only if $\Gamma$ contains a section $\Sigma$ such that $\Sigma=\bigoplus\limits_{M\in\Sigma}M$ is a $\tau$-rigid module and $|A|=|\Sigma|$.
\end{theorem}

\begin{proof}
Because of Lemma \ref{caminito}, it suffices to show that there is no path in $\Gamma$ from $\Sigma$ to a projective not in $\Sigma$, and no path from an injective not in $\Sigma$ to $\Sigma$.

Assume that there exists a path in $\Gamma$ from a module $Y\in \Sigma$ to a projective $P$:
$$Y=Y_0\rightarrow Y_1\rightarrow\dots\rightarrow Y_t=P$$
Because $\Sigma$ is sincere, there exists $X\in\Sigma$ such that Hom$_A(P,X)\neq 0$. Because $\Gamma$ is generalized standard, this implies the existence of a path in $\Gamma$
$$P=X_0\rightarrow X_1\rightarrow\dots\rightarrow X_s=X$$
Composing both paths we get a path in $\Gamma$ from $Y$ to $X$ passing through $P$. Convexity of $\Sigma$ in $\Gamma$ implies that $P\in\Sigma$. This proves the first statement. The second is dual.  

Conversely let $\Sigma$ be a complete slice in mod$A$. Then $\Sigma$ is a $\tau$-rigid section and $|\Sigma|=|A|$ because $\Sigma$ is a tilting module. Moreover $\Sigma$ lies in the connecting component of $\Gamma_A$, which is convex and generalized standard, see \cite[Chapter VIII]{AsSS}. This completes the proof. 
\end{proof}

Note that in the previous theorem, the local condition $Ann\Sigma=0$ in the criterion of Liu-Skowronski (see \cite[Chapter VIII, Section 5]{AsSS}) is replaced by a global one, namely the convexity of the connected component $\Gamma$.

\begin{prop}\label{proposicioncompleta}
Let $A$ be a tilted algebra. Then $\Sigma$ is a complete $\tau$-slice if and only if $\Sigma$ is a complete slice.
\end{prop}

\begin{proof}
It follows from Lemma \ref{tiltedcociente} that complete slices are precisely faithful $\tau$-slices. 

Conversely, let $\Sigma$ be a complete $\tau$-slice, $\tilde{\Sigma}$ a complete slice in mod$A$ and $(Fac\tilde{\Sigma}, Sub(\tau\tilde{\Sigma}))$ be the torsion pair associated to the latter. We consider three cases. 

First suppose that every $M\in\Sigma$ belong to $Sub(\tau\tilde{\Sigma})$. Then $pd\Sigma \leq 1$ by \cite[Lemma VIII.3.2]{AsSS}. Then $\Sigma$ is a tilting module by Proposition \ref{propix}. Hence $\Sigma$ is faithful.

Dually suppose that every $M\in\Sigma$ belong to $Fac\tilde{\Sigma}$. Then $id\Sigma \leq 1$. Hence $\Sigma$ is a cotilting module. Hence $\Sigma$ is faithful.

Finally, suppose that $\Sigma$ is neither contained in $Sub(\tau\tilde{\Sigma})$ nor in $Fac\tilde{\Sigma}$. Then there exists indecompossable direct summands $Y_i,Y_j$ of $\Sigma$ such that $Y_i\in Fac\tilde{\Sigma}$ and $Y_j\in Sub(\tau\tilde{\Sigma})$ because $\tilde{\Sigma}$ induces a splitting torsion pair in mod$A$ by \cite[Lemma.VIII.3.2]{AsSS}. Hence $\Sigma$ is contained in the connecting component $\Gamma$ containing $\tilde{\Sigma}$. Therefore $\Gamma$ is convex and generalized standard. It is known that $\Gamma$ is a full connected subquiver of $\mathbb{Z}\tilde{\Sigma}$. Thus $\Sigma$ is a presection in the stable quiver $\mathbb{Z}\tilde{\Sigma}$ because it is a presection in $\Gamma$. It is shown in \cite[Proposition 7]{ABS2} that presections and sections coincide in a stable translation quiver. Therefore $\Sigma$ is a section in $\mathbb{Z}\tilde{\Sigma}$. 

We claim that $\Sigma$ is a section in $\Gamma$. First, $\Sigma$ is convex and acylic in $\Gamma$ because it is convex and acyclic in $\mathbb{Z}\tilde{\Sigma}$. Because $\Sigma$ is a section in $\mathbb{Z}\tilde{\Sigma}$, then for each $X\in\Gamma\subseteq\mathbb{Z}\tilde{\Sigma}$ there exists a unique $Y_k\in\Sigma$ and a unique $n\in\mathbb{Z}$ such that $X=\tau^nY_k$. Therefore $\Sigma$ is a section in $\Gamma$. Moreover $\Sigma$ is $\tau$-rigid because $\Sigma$ is a $\tau$-tilting module and $|A|=|\Sigma|$. Then Theorem \ref{taurigidsection} implies that $\Sigma$ is a complete slice in mod$A$.
\end{proof}

\begin{cor}
Let $A$ be an algebra having a complete $\tau$-slice $\Sigma$ in $\Gamma_A$. Then $A$ is tilted having $\Sigma$ as a complete slice if and only if $\Sigma$ is faithful.
\end{cor}

\begin{proof}
If $\Sigma$ is faithful, then Corollary \ref{taurodajasfieles} implies that $A$ is tilted. Conversely, if $A$ is tilted then $\Sigma$ is a complete slice by Proposition \ref{proposicioncompleta}. So $\Sigma$ is faithful. 
\end{proof}

The following example shows that the hypothesis of $\Gamma$ being convex cannot be removed from Theorem \ref{taurigidsection}.

\begin{ex}
Consider the algebra $A$ given by the path algebra of the quiver 
$$
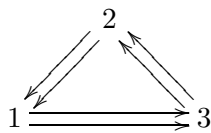
\xymatrix{
 & 2\ar@<0.5 ex>[dl]\ar@<-0.5 ex>[dl] & \\
 1\ar@<0.5 ex>[rr]\ar@<-0.5 ex>[rr] & & 3\ar@<0.5 ex>[ul]\ar@<-0.5 ex>[ul] }$$
factored by its square radical. Consider the connected component $\Gamma$ of the Auslander-Reiten quiver of $A$ represented in Figure 2. We can see that the complete $\tau$-slice $\Sigma=\rep{11\\2}\oplus\rep{1}\oplus\rep{3\\11}$ is a $\tau$-rigid section in $\Gamma$. It is also true that $\Gamma$ is generalized standard. However $A$ is not tilted and $\Sigma$ is not a complete slice.
\end{ex}

\begin{figure}[H]
$$\xymatrix{
 & & & & & \rep{3\\11}\ar@<0.5 ex>[dr]\ar@<-0.5 ex>[dr] & &\rep{333\\1111}\ar@<0.5 ex>[dr]\ar@<-0.5 ex>[dr]& \\
 \dots\ar@<0.5 ex>[dr]\ar@<-0.5 ex>[dr]& &\rep{111\\22}\ar@<0.5 ex>[dr]\ar@<-0.5 ex>[dr] & &\rep{1}\ar@<0.5 ex>[ur]\ar@<-0.5 ex>[ur] & &\rep{33\\111}\ar@<0.5 ex>[ur]\ar@<-0.5 ex>[ur] & & \dots \\
 &\rep{1111\\222}\ar@<0.5 ex>[ur]\ar@<-0.5 ex>[ur] & &\rep{11\\2}\ar@<0.5 ex>[ur]\ar@<-0.5 ex>[ur] & & & & &
 }$$
\caption{}
\end{figure}

Let $A$ be an algebra. For every presentation $A=kQ/I$, one compute its fundamental group $\pi_1(Q,I)$ (see \cite{M-VdlP}).  We say that a triangular algebra $A$ is \textit{simply connected} if the fundamental group $\pi(Q,I)$ is trivial for every presentation of $A=kQ/I$. We can also associate a fundamental group $\pi_1(\Gamma)$ to a connected component $\Gamma$ of the Auslander-Reiten quiver of $A$ (see \cite{BG}). In particular, a component $\Gamma$ is simply connected if and only if its orbit graph is a tree. Moreover, if a tilted algebra is simply connected then its connecting component is simply connected, see \cite{LM}. For the rest of the section we study simply connected algebras.

\begin{lem}\label{compsimple}
Let $\Gamma$ be a connected component of the Auslander-Reiten quiver $\Gamma_A$ of an algebra $A$. Assume that $\Gamma$ has a $\tau$-slice $\Sigma$. If $\Gamma$ is simply connected then each point of $\Sigma$ belongs to a different $\tau$-orbit of $\Gamma$.
\end{lem}

\begin{proof}
First of all, because $\Gamma$ is simply connected, $\Gamma$ is acyclic. 

Let $\Omega$ be a full connected subquiver of $\Sigma$ and consider a numeration in $\Omega_0=\{X_1, X_2, \dots, X_t\}$ such that there is an arrow from $X_i$ to $X_{i+1}$ or an arrow from $X_{i+1}$ to $X_{i}$ for every $i$. 

Suppose that $\Omega_0=\{X_1, X_2\}$. We can suppose without loss of generality that there is an arrow from $X_1$ to $X_2$. Suppose that $X_1$ and $X_2$ are in they same $\tau$-orbit. Then there exists $n \in \mathbb{N}$ such that $X_1 = \tau^n X_2$ or $X_1= \tau^{-n} X_2$. If $X_2= \tau^n X_1$, then there is a cycle in $\Gamma$, a contradiction to our hypothesis. If $X_2 = \tau^{-n} X_1$, then there is an arrow from $X_2$ to $\tau^{-1} X_1 = \tau^{n-1} X_2$. So we can apply the same argument to deduce that there is a cycle in $\Gamma$. 

Suppose that $\Omega_0=\{X_1, X_2, X_3\}$ and that $X_3 = \tau^n X_1$ with $n\geq 1$. Notice that the last assumption implies that $X_3$ is not injective and, dually, that $X_1$ is not projective. We have, without loss of generality, the following four cases:
\begin{enumerate}
    \item There is an arrow from $X_1$ to $X_2$ and one arrow from $X_2$ to $X_3$: In this case it is trivial that there is a cycle in the component;
    \item There is an arrow from $X_2$ to $X_1$ and one arrow from $X_2$ to $X_3$: In this case there exists an arrow from $\tau_A X_1$ to $X_2$ and we are back to the first case, and we can construct a cycle;
    \item There is an arrow from $X_1$ to $X_2$ and one arrow from $X_3$ to $X_2$: If we consider $X_1$, $X_2$ and $\tau_A^{-1} X_3$, we are in the situation of the first case once again, where a cycle in the component arises;
    \item There is an arrow from $X_3$ to $X_2$ and one arrow from $X_2$ to $X_1$: because $\Sigma$ is a presection we know that $X_3$ is not $\tau_A X_1$, then considering $\tau X_1$, $X_2$ and $X_3$ we have the situation of the third case and once again we can construct a cycle in $\Gamma$. 
\end{enumerate}
Therefore if $X_1$ and $X_3$ belong to the same $\tau$-orbit then there exists a cycle in $\Gamma$, a contradiction with our hypothesis. 

Finally, suppose that $\Omega_0=\{X_1, X_2, \dots, X_t\}$ with $t\geq 4$ and that $X_i$ does not belong to the $\tau$-orbit of $X_j$ if $1\leq|i-j|\leq t-2$. If $X_1$ belongs to the $\tau$-orbit of $X_t$ for some $\Omega$, then the orbit graph of $\Gamma$ contains a subgraph which is a cycle, contradicting the hypothesis of simple connectedness. Therefore every point of $\Sigma$ belongs to a different $\tau$-orbit.
\end{proof}

\begin{theorem}\label{simplyconnectedtilted}
Let $A$ be an algebra and $\Gamma_A$ be its Auslander-Reiten quiver. Assume that $\Gamma_A$ has a connected component $\Gamma$ which is simply connected, convex and generalized standard. Then $\Sigma$ is a complete $\tau$-slice in $\Gamma$ if and only if $A$ is a simply connected tilted algebra having $\Sigma$ as a complete slice in mod$A$.
\end{theorem}

\begin{proof}
Necesity. Following Lemma \ref{caminito}, it is enough to show that there is no path in $\Gamma$ from $\Sigma$ to a projective and no path from an injective to $\Sigma$.

Fix a numbering in the points of $\Sigma_0=\{X_1, \dots, X_n\}$. Suppose that there exists a path $\omega'$ 
$$X\rightarrow Y_1 \rightarrow \dots\rightarrow Y_t\rightarrow P $$
in $\Gamma$ from $X\in \Sigma$ to an indecomposable projective module $P$ which does not belong to $\Sigma$. Because $\Sigma$ is a finite presection we can suppose without loss of generality the existence of $X_1\in\Sigma$ such that $Y_1=\tau^{-1}X_1$.

The $\tau$-slice $\Sigma$ is a sincere $A$-module because it is a complete $\tau$-slice. Then there exists $Y \in \Sigma_0$ such that Hom$_A(P,Y)\neq 0$. This implies the existence of a path in $\Gamma$ from $P$ to $Y$ because $\Gamma$ is generalized standard. Let $\omega''$ be such a path:
$$P \rightarrow Y_{t+1}\rightarrow \dots \rightarrow Y_{m-1} \rightarrow  Y_m \rightarrow Y.$$
Because $\Sigma$ is a finite presection, \cite[Corollary IX.2.3]{AsSS} implies that there exist $j$ and $X_k\in\Sigma$ such that $t+1\leq j\leq m$ and $Y_j=\tau X_k$. We can suppose without loss of generality that $j=m$.  

Take the composition path $\omega=\omega'\omega''$.
$$X\rightarrow Y_1 \rightarrow \dots \rightarrow P \rightarrow \dots \rightarrow Y_{m-1} \rightarrow  Y_m \rightarrow Y$$
Consider the set $\mathcal{I}=\{Y_i = \tau^{z_i} X_h : Y_i \in \omega, X_h \in \Sigma, z_i \in \mathbb{Z} \}$. This set is not empty because $Y_1, Y_m \in \mathcal{I}$. Also, if, for each $h$ with $1\leq h\leq n$ we let $\mathcal{I}_h=\{Y_i = \tau^{z_i} X_h : Y_i \in \omega, z_i \in \mathbb{Z} \}$, then
$$\mathcal{I}=\bigcup_{h=1}^{n} \mathcal{I}_h$$
Moreover, Lemma \ref{compsimple} implies that $\mathcal{I}_i\cap\mathcal{I}_j=\emptyset$.

Let $Y_s$ be the first point in the path $\omega$ such that $Y_s=\tau^{z_s} X_s$ with $z_s\geq 1$ for some $X_s\in\Sigma_0$. Such a $Y_s$ exists because $Y_m=\tau X_k$. If $Y_s=Y_1$ then $Y_s=\tau^{z_s} X_1$ by Lemma \ref{compsimple}. This allows us to construct a cycle from $X_1$ to itself passing by $\tau^{-1} X_1$ and $\tau^{z_s} X_1$,
$$X_1\rightarrow *\rightarrow \tau^{-1}X_1=\tau^{z_s}X_1\rightarrow *\rightarrow \dots\rightarrow \tau X_1\rightarrow *\rightarrow X_1$$
a contradiction to our hypothesis.

Suppose that $s>1$ and take $Y_r=\tau^{z_r}X_r$ such that $r<s$ and $Y_i\not\in\mathcal{I}$ if $r<i<s$. Then there exists a subpath $\tilde{\omega}$ of $\omega$ from $Y_r$ to $Y_s$. Because $\Gamma$ is simply connected, then either $X_r=X_s$ or they lie in two neighboring $\tau$-orbtis. If $X_r=X_s$ then it is possible to construct a cycle from $X_r$ to itself as follows. 

$$X_r\rightarrow *\rightarrow \tau^{-1}X_r\rightarrow\dots\rightarrow\tau^{z_r}X_r\overset{\tilde{\omega}}\rightsquigarrow\tau^{z_s}X_r \rightarrow \dots\rightarrow \tau X_r\rightarrow *\rightarrow X_r$$

Otherwise $X_r$ and $X_s$ belong to neighboring $\tau$-orbits because $\Gamma$ is simply connected. Moreover, because $\Sigma$ is connected then there is an arrow from $X_r$ to $X_s$ (or from $X_s$ to $X_r$). Therefore we can construct a cycle from $X_r$ to itself passing by $Y_r$, $Y_s$ and $\tau X_s$ (or $X_s$),
$$X_r\rightarrow *\rightarrow \tau^{-1}X_r\rightarrow\dots\rightarrow\tau^{z_r}X_r\overset{\tilde{\omega}}\rightsquigarrow\tau^{z_s}X_s \rightarrow \dots\rightarrow \tau X_s\rightarrow X_r$$
a contradiction to our assumption. 

Hence there is no path from $\Sigma$ to a projective $P$ in $\Gamma$ which is not in $\Sigma$. Dually, there is no path from an injective $I$ in $\Gamma$ not in $\Sigma$ to $\Sigma$. Therefore Lemma \ref{caminito} implies that $A$ is tilted and $\Sigma$ is a complete slice in mod$A$. 

Moreover $\Sigma$ is a tree because $\Gamma$ is simply connected, see \cite{BG}. Then \cite[Theorem A]{LM} implies that $A$ is simply connected.

Sufficiency. Let $\Sigma$ be a complete slice in mod$A$. Then $\Sigma$ is a complete $\tau$-slice lying in the connecting component $\Gamma$ of $\Gamma_A$, which is convex and generalized standard. Moreover $\Gamma$ is simply connected by \cite[Theorem 4.3]{M-VdlP}. This finishes the proof.
\end{proof}

The following example shows that the hypothesis of $\Gamma$ being simply connected cannot be removed in Theorem \ref{simplyconnectedtilted}.

\begin{ex}
Consider the algebra $A$ being the path algebra of the quiver 
$$
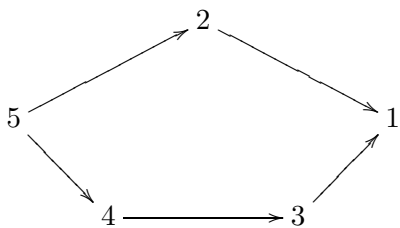
\xymatrix{
 & & 2\ar[rrd] & & & \\
 5\ar[rru]\ar[rd]& & & &1 \\
 &4\ar[rr] & &3\ar[ru] & }$$
modulo its square radical. See its Auslander-Reiten quiver in Figure 3. We can see that $\Gamma_{A}$ has only one connected component which is acyclic and trivially generalized standard but not simply connected. Here $\Sigma=\rep{2\\1}\oplus\rep{23\\1}\oplus\rep{2}\oplus\rep{5\\42}\oplus\rep{5\\2}$ is a complete $\tau$-slice. However $A$ is not a tilted algebra and $\Sigma$ is not a complete slice. Note that $\rep{2\\1}$ and $\rep{5\\2}$ belong to the same $\tau$-orbit.
\end{ex}

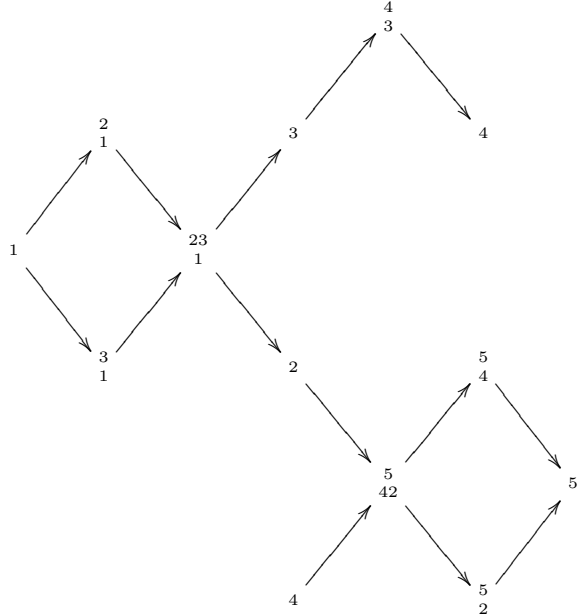
\begin{figure}\label{ARQA4}
$$
\xymatrix{
 & & & & \rep{4\\3}\ar[rd] & & \\
 & \rep{2\\1}\ar[rd] & & \rep{3}\ar[ru] & & \rep{4} & \\ 
  \rep{1}\ar[rd]\ar[ru]& & \rep{23\\1}\ar[ru]\ar[rd] & & & & \\ 
 & \rep{3\\1}\ar[ru] & & \rep{2}\ar[rd] & & \rep{5\\4}\ar[rd] & \\ 
 & & & & \rep{5\\42}\ar[rd]\ar[ru] & & \rep{5} \\ 
 & & & \rep{4}\ar[ru] & & \rep{5\\2}\ar[ru] &  
}$$
\caption{The Auslander-Reiten quiver of $A$}
\end{figure}

\section{Some Extensions of Algebras with $\tau$-Slices}

In this final section we take an algebra $A$ with a $\tau$-slice $\Sigma$ and we construct some new algebras $\tilde{A}$ based on $A$ and having a $\tau$-slice in mod$\tilde{A}$. First we consider one-point (co)extensions. Afterwards we consider split-by-nilpotent extensions. Finally we give a series of algebras having $\tau$-slices in order to illustrate the results of this section and Theorem \ref{teorema1}.

Let $A$ be an algebra and $M$ an $A$-module. We denote by $B=A[M]$ the \textit{one-point extension} of $A$ by $M$, that is, the matrix algebra 
$$B=\left(\begin{matrix} 
A & 0 \\
M & k
\end{matrix}\right)$$
with the ordinary matrix addition and the multiplication induced from the module structure of $M$. Note that the ordinary quiver of the algebra $B$ is obtained from that of $A$ by adding a new vertex, called the \textit{extension point}, such that the radical of the indecomposable projective module associated to this vertex is $M$. 

Following \cite{AZ1} we denote by $\mathcal{C}_A$ the class of objects $X$ in mod$A$ such that the almost split sequences ending with $X$ in mod$A$ and in mod$B$ coincide. 

\begin{prop}\label{onepoint}
Let $A$ be an algebra with a complete $\tau$-slice $\Sigma$ and let $X \in add\Sigma$. Then the algebra $B=A[X]$ has a complete $\tau$-slice $\tilde{\Sigma}=\Sigma \oplus P_x$, where $P_x$ is the projective $B$-module associated to the extension point.
\end{prop} 

\begin{proof}
Set $\tilde{\Sigma} = \Sigma \oplus P_x$.

First note that every indecomposable projective $A$-module is also an indecomposable projective $B$-module. Moreover we have Hom$_A(X,\tau_A \Sigma)=0$, because $X$ belongs to add$(\Sigma)$. Therefore every indecomposable non-projective summand of $\Sigma$ belongs to $\mathcal{C}_A$ by \cite[Corollary 2.6]{AZ1}. Consequently $\tau_A \Sigma \cong \tau_B \Sigma$.

We have Hom$_B(\tilde{\Sigma},\tau_B\tilde{\Sigma})=$Hom$_B(P_x,\tau_B \Sigma) \oplus$ Hom$_B(\Sigma, \tau_B \Sigma)$ by the additivity of the Hom functor and the fact that $P_x$ is a projective $B$-module. Because $\tau_B \Sigma$ is an $A$-module, the simple module $S_x$ associated to the extension point $x$ is not a composition factor of $\tau_B \Sigma$. Then Hom$_B(P_x, \tau_B \Sigma)=0$. Also, Hom$_B(\Sigma, \tau_B \Sigma)=$Hom$_A(\Sigma, \tau_A \Sigma)=0$. Hence $\tilde{\Sigma}$ is a $\tau$-rigid $B$-module.

Now we prove that $\tilde{\Sigma}$ is a presection. Let $f: M \rightarrow N$ be an arrow in $\Gamma_B$ and suppose that $N$ belongs to $\tilde{\Sigma}$. If $N=P_x$ then $M \in \tilde{\Sigma}$ because the minimal right almost split morphism ending in a projective is the inclusion of the radical and by hypothesis every direct summand of rad$P_x$ belong to add$\Sigma$. If $N \in \Sigma$ then $f$ is an arrow in $\Gamma_A$ because the almost split sequences ending with $N$ in mod$A$ and in mod$B$ coincide. We know that $\Sigma$ is a presection in mod$A$. Then either $M \in \Sigma$ or $M'=\tau^{-1}_A M \in \Sigma$. Therefore either $M$ or $\tau_B^{-1} M$ belong to $\tilde{\Sigma}$. 

Now, suppose that $M \in \tilde{\Sigma}$ and that $N$ is a projective $B$-module. If $N=P_x$ then $N \in \tilde{\Sigma}$. If $N$ is another projective $B$-module, then $N$ is a projective $A$-module. By Lemma \ref{presecpro}, $N \in \Sigma$. Hence $N \in \tilde{\Sigma}$. 

Assume now that $N$ is not projective. Then $\tau_B N \neq 0$. Therefore there is an arrow $f': \tau_B N \longrightarrow M$. So, either $\tau_B N$ belongs to $\Sigma$ or $\tau^{-1}_A \tau_B N$ belongs to $\Sigma$. If $\tau_B N$ belongs to $\Sigma$, then there is nothing to prove. If $\tau^{-1}_A \tau_B N$ is in $\Sigma$, then in particular $\tau_B N$ is an $A$-module. Therefore $\tau_A N =\tau_B N$. So $N = \tau^{-1}_A \tau_A N = \tau^{-1}_A \tau_B N$ belongs to $\Sigma$.

Then $\tilde{\Sigma}$ is a $\tau$-rigid presection in mod$B$. Therefore Lemma \ref{rigidezalcanza} implies that $\tilde{\Sigma}$ is a $\tau$-slice in mod$B$. Moreover it is a complete $\tau$-slice because $|\tilde{\Sigma}|=|\Sigma|+1=|A|+1=|B|$.
\end{proof}

As a corollary we get the following well-known result, see for instance \cite[Lemma 3.5]{OS}.

\begin{cor}
Let $A$ be a tilted algebra with complete slice $\Sigma$ and $X \in add\Sigma$. If we denote by $P_x$ the projective $B$-module associated to the extension point, then the algebra $B=A[X]$ is a tilted algebra having $\tilde{\Sigma}=\Sigma \oplus P_x$ as a complete slice.
\end{cor}

\begin{rmk}
Similar results have been proved by Oryu and Schiffler in \cite{OS} working in the context of cluster tilted algebras. In particular they proved that the one-point extension $B[P]$ of a cluster tilted algebra $B$ by a projective module $P$ lying in a local slice is cluster tilted (see \cite[Theorem 3.6]{OS}).
\end{rmk}

\begin{lem}\label{propita}
Let $\Sigma$ be a $\tau$-slice in mod$A$. Then Fac$(\tau^{-1}\Sigma)\cup$add$\Sigma=$Fac$(\Sigma)$.
\end{lem}

\begin{proof}
In Corollary \ref{tauiguales} we showed that $\tau^{-1}_A \Sigma= \tau^{-1}_{C}\Sigma$, where $C=A/Ann\Sigma$. Moreover $\Sigma$ is a complete slice in mod$C$ and $C$ is a tilted algebra by Corollary \ref{tiltedcociente}. Therefore it is sufficient to show the statement in the case where $A$ is tilted.

Let $A$ be tilted and $\Sigma$ be a complete slice in mod$A$. Then \cite[Lemma VIII.3.2]{AsSS} implies that $\tau^{-1}_A \Sigma\in$ Fac$\Sigma$ because $\tau^{-1}_A\Sigma$ is a successor of $\Sigma$. Therefore Fac$(\tau^{-1}\Sigma)\subseteq$Fac$(\Sigma)$. Moreover, for every indecomposable module $N\in$ Fac$\Sigma$ there exists a epimorphism $\pi: M\rightarrow N$, where $M\in$ add$\Sigma$. If $N$ is not a direct summand of $\Sigma$, $\pi$ factors through $\tau^{-1}_A \Sigma$ by \cite[Lemma VIII.5.4]{AsSS}. Therefore $N\in$ Fac$(\tau^{-1}_A\Sigma)$.
\end{proof}

\begin{cor}\label{corex}
Let $A$ be an algebra with a $\tau$-slice $\Sigma$ and let $X\in Fac(\tau^{-1}_A\Sigma)$. Then $\Sigma$ is a $\tau$-slice for the algebra $B=A[X]$.
\end{cor}

\begin{proof}
By Proposition \ref{propita}, $Fac(\tau^{-1}_A \Sigma)$ is contained in $Fac_A \Sigma$. Therefore Hom$_A(X,\tau_A \Sigma)=0$. Then \cite[Corollary 2.6]{AZ1} implies that $Y\in \mathcal{C}_A$ for every $Y$ indecomposable direct summand of $\Sigma$. Hence Hom$_B(\Sigma,\tau_B \Sigma)=$Hom$_A(\Sigma, \tau_A \Sigma)=0$, in other words, $\Sigma$ is a $\tau$-rigid $B$-module. Also Hom$_A(X,\Sigma)=0$ by Corollary \ref{tauiguales}. Hence $Z\in \mathcal{C}_A$ for every indecomposable direct summand $Z$ of $\tau^{-1}_A\Sigma$. Then every almost split sequence starting or ending in $\Sigma$ in mod$A$ coincides with the corresponding almost split sequence starting or ending in $\Sigma$ in mod$B$. Therefore $\Sigma$ is a presection in $\Gamma_B$ because it is a presection in $\Gamma_A$.

The assertion follows from Lemma \ref{rigidezalcanza}.
\end{proof}

Using the dual arguments we can prove the following results. We denote by $[X]A$ the one point coextension of an algebra $A$ by the $A$-module $X$.

\begin{theorem}
Let $A$ be an algebra with a complete $\tau$-slice $\Sigma$ and let $X \in add\Sigma$. Then the algebra $B=[X]A$ has a complete $\tau$-slice $\tilde{\Sigma}=\Sigma \oplus I_x$, where $I_x$ is the injective $B$-module associated to the extension point.
\end{theorem} 

\begin{cor}\label{corcoex}
Let $A$ be an algebra with a $\tau$-slice $\Sigma$ and let $X \in Sub(\tau_A \Sigma)$. Then $\Sigma$ is a $\tau$-slice for the algebra $B=[X]A$.
\end{cor}

\begin{rmk}
Note that in Corollary \ref{corex} and in Corollary \ref{corcoex} the $\tau$-slice $\Sigma$ will never be complete because $|\Sigma|=|A|$ and $|B|=|A| + 1$.
\end{rmk}

Let $A$ be an algebra and $Q$ an $A$-$A$-bimodule equipped an $A$-$A$-morphism $\mu:Q\otimes_A Q\rightarrow Q$. Then $B=A\oplus Q$ becomes an algebra if one defines the multiplication by $$(a,q).(a',q')=(aa',aq'+qa'+qq')$$ Moreover there exists a split exact sequence of $k$-vector spaces
$$0\rightarrow Q\overset{\iota}\rightarrow B\overset{\pi}\rightarrow A\rightarrow 0$$
where $\iota:q\mapsto (0,q)$ is the inclusion of $Q$ as a two-sided ideal of $B=A\oplus Q$, and the projection (algebra) morphism $\pi:(a,q)\mapsto a$ has as section the inclusion (algebra) morphism $\sigma: a\mapsto (a,0)$. Then we say that $B$ is the \textit{split extension} of $A$ by $Q$. If $Q$ is nilpotent we say that $B$ is the \textit{split-by-nilpotent extension} of $A$ by $Q$. In particular, if $\mu(q\otimes q')=0$ for every $q,q'\in Q$, we say that $B$ is the \textit{trivial extension} of $A$ by $Q$. We refer the reader to \cite{AZ1} for further details.

For the rest of the paper we will consider split-by-nilpotent extensions.

The following theorem gives a necessary and sufficient condition on $Q$ to make $\Sigma$ a complete $\tau$-slice in mod$B$.

\begin{theorem}\label{splitex}
Let $\Sigma$ be a complete $\tau$-slice in mod$A$ and $Q$ a nilpotent $A$-$A$-bimodule. Consider the split-by-nilpotent extension $B$ of $A$ by $Q$. 
Then $\Sigma$ is a complete $\tau$-slice in mod$B$ if and only if $Q_A\in Fac(\tau^{-1}_A \Sigma)$ and $D(_A Q)\in$ Sub$(\tau_A \Sigma)$. Moreover, if $A$ is a tilted algebra and $\Sigma$ is a complete slice in mod$A$ then Ann$_B \Sigma=Q$.
\end{theorem}

\begin{proof}
Let $M$ be an indecomposable direct summand of $\Sigma$. Suppose that $Q_A\in Fac(\tau^{-1}_A \Sigma)$ and $D(_A Q)\in$Sub$(\tau_A \Sigma)$. Then Hom$_A(M,DQ) =0$ and Hom$_A(Q,\tau_AM) =0$. If $M$ is an indecomposable projective $A$-module, then $M=eA$ for some primitive idempotent $e$ of $A$. Therefore the simple module top$(eA)$ is not a composition factor of $D(_A Q)$. Now $D(eQ)=D(Q)e=0$. Hence $eQ=0$. Consequently \cite[Corollary 1.4]{AZ1} implies that $M$ is a projective $B$-module. If $M$ is not a projective $A$-module \cite[Theorem 2.1]{AZ1} implies that the almost split sequences ending with $M$ in mod$A$ and in mod$B$ coincide. Hence the almost split sequences ending with $M$ in mod$A$ and in mod$B$ coincide for all modules $M \in \Sigma$.   

Dually we can prove that every indecomposable injective summand of $\Sigma$ is also an injective $B$-module. Also the dual statement of \cite[Theorem 2.1]{AZ1} implies that the almost split sequences starting with $M$ in mod$A$ and in mod$B$ coincide for all module $M \in \Sigma$ under the hypothesis that $Q_A$ belongs to Fac$_A(\tau^{-1}_A \Sigma)$ and and $D(_A Q)$ belongs to $Sub(\tau_A \Sigma)$ because Hom$_A(\tau^{-1}_AM,DQ) =0$ and Hom$_A(Q,M) =0$. Therefore $\Sigma$ is a presection in $\Gamma_B$ because it is a presection in $\Gamma_A$. 
Also, Hom$_B(\Sigma, \tau_B \Sigma)=$Hom$_A(\Sigma,\tau_A \Sigma)=0$. Thus $\Sigma$ is a $\tau$-slice in $\Gamma_A$, which is complete because $|A|=|B|$. 

Conversely, suppose that $\Sigma$ is a complete $\tau$-slice in mod$B$. First we prove that $Q\subseteq Ann\Sigma$. We have that $B \cong A \oplus Q$ as $k$-vector spaces. The action of $b=(a,q)\in B$ on an $A$-module $M$ is given by the following formula.
$$M b = M (a,q) := M a$$
Then $Q=\{(0,q)\}\subseteq Ann_BM$. Therefore Corollary \ref{tauiguales} implies that $\tau_B \Sigma = \tau_A \Sigma$. Hence $Hom_C(M, D_CQ)=0$ for every $M$ non-projective indecomposable summand of $\Sigma$ by \cite[Theorem 2.1]{AZ1}. If $M=eA$ is a projective summand of $\Sigma$, then Theorem \ref{teorema1} implies that $M=eA=eC$. Moreover $D(Qe)=0$ by \cite[Corollary 1.4]{AZ1}. Then Hom$_A(eA,DQ)=0$. Therefore Hom$_A(\Sigma, DQ)=0$. Thus $DQ \in Sub_A(\tau_A \Sigma)$. The other assertion is proved dually. 

Finally, if $A$ is tilted and $\Sigma$ is a complete slice, then $\Sigma$ is a faithful $A$-module by Corollary \ref{taurodajasfieles} and Ann$_B \Sigma = Q$, finishing the proof.
\end{proof}

In \cite{ABS1}, Assem, Br\"ustle and Schiffler showed that a given cluster tilted algebra $\tilde{C}$ is the trivial extension of a tilted algebra $C$ by the $C$-$C$-bimodule $E=Ext_C^2(DC,C)$. They called the extension by this particular bimodule the \textit{relation extension}. In the subsequent paper \cite{ABS2}, they showed that every complete slice $\Sigma$ in mod$C$ become a local slice in mod$\tilde{C}$ and, conversely, that every local slice in mod$\tilde{C}$ arises in this way. Later on, Assem, Bustamante, Dionne, Le Meur and Smith considered in \cite{ABDLS} \textit{partial relation extensions}. A partial relation extension of a tilted algebra $C$ is the trivial extension of $C$ by a direct summand $E'$ of $E$ as $C$-$C$-bimodule. They proved that every partial relation extension of a tilted algebra $C$ has a local slice (see \cite[Theorem 3]{ABDLS}). As a consequence of the previous theorem and the explicit computation of the $C$-module structure $E$ and $D(E)$, made by Schiffler and Serhiyenko in \cite{SS}, we get the following corollary, implying \cite[Theorem 19]{ABS2} and \cite[Theorem 3]{ABDLS}.

\begin{cor}\label{corsplit}
Let $C$ be a tilted algebra, $\Sigma$ a complete slice in mod$C$, $E$ as above and $\tilde{E}\in$add$_{C-C}E$. Then $\Sigma$ is a complete $\tau$-slice in mod$\tilde{C}$, where $\tilde{C}$ is the trivial extension of $C$ by $\tilde{E}$.
\end{cor}

\begin{proof}
By Theorem \ref{splitex}, it is enough to prove that $\tilde{E}_C\in$Fac$(\tau^{-1}\Sigma)$ and $D(\tilde{E})_C\in$Sub$(\tau\Sigma)$. We prove the first statement and the second is dual. 

It was shown in \cite[Proposition 4.1]{SS} that $E_C\cong \tau^{-1}\Omega^{-1} C_C$, where $\Omega^{-1}$ stands for the cosyzygy functor. Let $N$ be a non-injective indecomposable direct summand of $\Omega^{-1}C$. Then there exists an indecomposable injective module $I_N$ such that Hom$_C(I_N,N)\neq 0$. Moreover, there exists an indecomposable $L\in\Sigma$ such that Hom$_C(L,I_N)\neq 0$ because $\Sigma$ is sincere. Therefore there is a path $L\rightarrow I_N\rightarrow N\rightarrow *\rightarrow\tau^{-1}N$ in mod$A$. Hence $\tau^{-1}N$ is a proper successor of $\Sigma$. Thus $N\in$ Fac$(\tau^{-1}\Sigma)$, implying that $E_C\in$ Fac$(\tau^{-1}\Sigma)$. This finishes the proof.
\end{proof}

Now we show algebras with $\tau$-slices. 
\begin{ex}

Consider the cluster tilted algebra $\tilde{A}$ given by the quiver 
$$
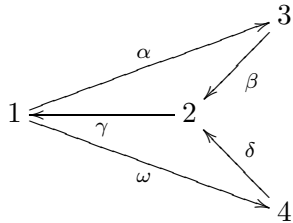
\xymatrix{
 & & & 3\ar[dl]^{\beta} \\
1\ar[urrr]^{\alpha}\ar[drrr]_{\omega} & & 2\ar[ll]^{\gamma} & \\
 & & &4\ar[ul]_{\delta} }$$
modulo the ideal $I=\langle \alpha\beta-\omega\delta, \beta\gamma, \delta\gamma, \gamma\omega, \gamma\alpha \rangle$. As we can see in Figure \ref{tAARQ}, the module $\Sigma=\rep{2\\1}\oplus\rep{2}\oplus\rep{3\\2}\oplus\rep{4\\2}$ is a complete $\tau$-slice in mod$\tilde{A}$. Its annihilator is $Ann_{\tilde{A}}\Sigma=\langle \alpha, \omega\rangle$. Take $I'=\langle \omega\rangle$, contained in $Ann_{\tilde{A}}\Sigma$, and consider the algebra $A=\tilde{A}/I'$. The algebra $A$ is the path algebra of the quiver
$$
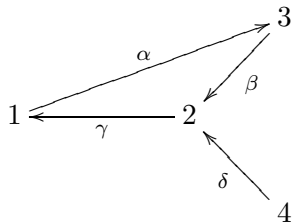
\xymatrix{
 & & & 3\ar[dl]^{\beta} \\
1\ar[urrr]^{\alpha} & & 2\ar[ll]^{\gamma} & \\
 & & &4\ar[ul]^{\delta} }$$
with radical square zero. Then, as proved in Theorem \ref{teorema1}, $\Sigma=\rep{2\\1}\oplus\rep{2}\oplus\rep{3\\2}\oplus\rep{4\\2}$ is a $\tau$-slice in mod$A$. See Figure \ref{AARQ}.

One can see that $A=A'[\rep{2}]$ where is $A'$ given by the quiver
$$
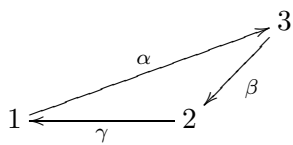
\xymatrix{
 & & & 3\ar[dl]^{\beta} \\
1\ar[urrr]^{\alpha} & & 2\ar[ll]^{\gamma} & }$$
with radical square zero and $\rep{2}$ is the simple module associated to the point $2$. See the Auslander-Reiten quiver of $A'$ in Figure 6. Then $\rep{2}$ belongs to the complete $\tau$-slice $\Sigma_1=\rep{2\\1}\oplus\rep{2}\oplus\rep{3\\2}$ and $\Sigma=\Sigma_1\oplus\rep{4\\2}$, where $\rep{4\\2}$ is the projective module associated to the extension point, agreeing with Theorem \ref{onepoint}. On the other hand we have that $\rep{2}\in Fac_{A'}(\tau^{-1}_{A'}\Sigma_2)$, where $\Sigma_2=\rep{1\\3}\oplus\rep{1}\oplus\rep{2\\1}$. Note that $\Sigma_2$ is a complete $\tau$-slice in mod$A'$ and a $\tau$-slice in mod$A$ which is not complete.  

Finally consider the algebra $C=A/\langle \alpha\rangle$. Corollary \ref{tiltedcociente} implies that $C$ is tilted and $\Sigma$ is a complete slice in mod$C$ because $Ann_A\Sigma=\langle\alpha\rangle$. The nilpotent ideal $Ann_A\Sigma$ of $A$ has an induced structure of $C-C$-bimodule. Moreover $A$ is the split-by-nilpotent extension of $C$ by $Ann_A\Sigma$. Also, as a right $C$-module $(Ann_A\Sigma)_C=\rep{3} \in Fac_C(\tau_C^{-1}\Sigma)$ and $D(_C(Ann_A\Sigma))=\rep{1}\in Sub\tau\Sigma$. Then Theorem \ref{splitex} implies that $\Sigma$ is a complete $\tau$-slice in mod$B$ induced by the complete slice $\Sigma$ in mod$C$.

\begin{figure}
    $$\xymatrix{ & & &\rep{3\\2}\ar[dr] & &\rep{4}\ar[dr] & &\rep{1\\3}\ar[dr] & \\
                \rep{1}\ar[r] &\rep{2\\1}\ar[r] &\rep{2}\ar[dr]\ar[ur] & &\rep{34\\2}\ar[dr]\ar[r]\ar[ur] &\rep{1\\34\\2}\ar[r] &\rep{1\\34}\ar[ur]\ar[dr] & &\rep{1} \\
                 & & &\rep{4\\2}\ar[ur] & &\rep{3}\ar[ur] & &\rep{1\\4}\ar[ur] &  }$$
    \caption{The Auslander-Reiten quiver of $\tilde{A}$}
    \label{tAARQ}
\end{figure}
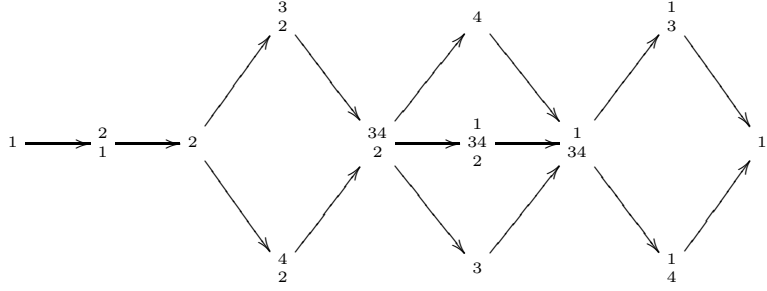

\begin{figure}
    $$\xymatrix{ &\rep{1\\3}\ar[dr] & & & &\rep{3\\2}\ar[dr]& &\rep{4} & \\
                \rep{3}\ar[ur] & &\rep{1}\ar[dr] & &\rep{2}\ar[r]\ar[ur] &\rep{4\\2}\ar[r] &\rep{34\\2}\ar[ur]\ar[r] &\rep{3}\ar[dr] & \\
                 & & &\rep{2\\1}\ar[ur] & & & & &\rep{1\\3} }$$
    \caption{The Auslander-Reiten quiver of $A$}
    \label{AARQ}
\end{figure}

\begin{figure}
    $$\xymatrix{ & &\rep{1\\3}\ar[dr] & & & &\rep{3\\2} \\
                 &\rep{3}\ar[ur] & &\rep{1}\ar[dr] & &\rep{2}\ar[ur] & \\
                \rep{3\\2}\ar[ur] & & & &\rep{2\\1}\ar[ur] & & }$$
    \label{A'ARQ}
    \caption{The Auslander-Reiten quiver of $A'$}
\end{figure}

\end{ex}

\newpage
\begin{center}
  \textsc{Acknowledgements}
\end{center}
The author thankfully acknowledges the financial support from the Facult\'e de Sciences of the Universit\'e de Sherbrooke. The author also thanks to Ibrahim Assem and David Smith for the attentive reading of the manuscript and their meaningful comments. Finally he is also grateful to Thomas Br\"ustle for useful discussions.

\def\cprime{$'$} \def\cprime{$'$}


\begin{thebibliography}{10}

\bibitem{AIR}
T.~Adachi, O.~Iyama, and I.~Reiten.
\newblock {$\tau$}-tilting theory.
\newblock {\em Compos. Math.}, 150(3):415--452, 2014.

\bibitem{ABS2}
I.~Assem, T.~Br{\"u}stle, and R.~Schiffler.
\newblock Cluster-tilted algebras and slices.
\newblock {\em J. Algebra}, 319(8):3464--3479, 2008.

\bibitem{ABS1}
I.~Assem, T.~Br{\"u}stle, and R.~Schiffler.
\newblock Cluster-tilted algebras as trivial extensions.
\newblock {\em Bull. Lond. Math. Soc.}, 40(1):151--162, 2008.

\bibitem{ABDLS}
I.~Assem, J.~C. Bustamante, J.~Dionne, P.~Le Meur, and D.~Smith.
\newblock Representation theory of partial relation extensions.
\newblock {\em arXiv}, math.RT(1604.01269), 2016.

\bibitem{AsSS}
I.~Assem, D.~Simson, and A.~Skowro{\'n}ski.
\newblock {\em Elements of the representation theory of associative algebras.
  {V}ol. 1}, volume~65 of {\em London Mathematical Society Student Texts}.
\newblock Cambridge University Press, Cambridge, 2006.
\newblock Techniques of representation theory.

\bibitem{AZ1}
I.~Assem and D.~Zacharia.
\newblock Full embeddings of almost split sequences over split-by-nilpotent
  extensions.
\newblock {\em Colloq. Math.}, 81(1):21--31, 1999.

\bibitem{APR}
M.~Auslander, M.~I. Platzeck, and I.~Reiten.
\newblock Coxeter functors without diagrams.
\newblock {\em Trans. Amer. Math. Soc.}, 250:1--46, 1979.

\bibitem{ARS}
M.~Auslander, I.~Reiten, and S.~O. Smal{\o}.
\newblock {\em Representation theory of {A}rtin algebras}, volume~36 of {\em
  Cambridge Studies in Advanced Mathematics}.
\newblock Cambridge University Press, Cambridge, 1997.
\newblock Corrected reprint of the 1995 original.

\bibitem{BGP}
I.~N. Bern{\v{s}}te{\u\i}n, I.~M. Gel{\cprime}fand, and V.~A. Ponomarev.
\newblock Coxeter functors, and {G}abriel's theorem.
\newblock {\em Uspehi Mat. Nauk}, 28(2(170)):19--33, 1973.

\bibitem{BG}
K.~Bongartz and P.~Gabriel.
\newblock Covering spaces in representation-theory.
\newblock {\em Invent. Math.}, 65(3):331--378, 1981/82.

\bibitem{BB}
S.~Brenner and M.~C.~R. Butler.
\newblock Generalizations of the {B}ernstein-{G}el\cprime fand-{P}onomarev
  reflection functors.
\newblock In {\em Representation theory, {II} ({P}roc. {S}econd {I}nternat.
  {C}onf., {C}arleton {U}niv., {O}ttawa, {O}nt., 1979)}, volume 832 of {\em
  Lecture Notes in Math.}, pages 103--169. Springer, Berlin-New York, 1980.

\bibitem{BMR}
A.~B. Buan, R.~J. Marsh, and I.~Reiten.
\newblock Cluster-tilted algebras.
\newblock {\em Trans. Amer. Math. Soc.}, 359(1):323--332 (electronic), 2007.

\bibitem{FZ1}
S.~Fomin and A.~Zelevinsky.
\newblock Cluster algebras. {I}. {F}oundations.
\newblock {\em J. Amer. Math. Soc.}, 15(2):497--529 (electronic), 2002.

\bibitem{HR}
D.~Happel and C.~M. Ringel.
\newblock Tilted algebras.
\newblock {\em Trans. Amer. Math. Soc.}, 274(2):399--443, 1982.

\bibitem{J}
G.~Jasso.
\newblock Reduction of {$\tau$}-tilting modules and torsion pairs.
\newblock {\em Int. Math. Res. Not. IMRN}, (16):7190--7237, 2015.

\bibitem{LM}
P.~Le~Meur.
\newblock Topological invariants of piecewise hereditary algebras.
\newblock {\em Trans. Amer. Math. Soc.}, 363(4):2143--2170, 2011.

\bibitem{LiuCrit}
S.~Liu.
\newblock Semi-stable components of an {A}uslander-{R}eiten quiver.
\newblock {\em J. London Math. Soc. (2)}, 47(3):405--416, 1993.

\bibitem{Liu}
S.~Liu.
\newblock Another characterization of tilted algebras.
\newblock {\em Arch. Math. (Basel)}, 104(2):111--123, 2015.

\bibitem{M-VdlP}
R.~Mart{\'{\i}}nez-Villa and J.~A. de~la Pe{\~n}a.
\newblock The universal cover of a quiver with relations.
\newblock {\em J. Pure Appl. Algebra}, 30(3):277--292, 1983.

\bibitem{OS}
M.~Oryu and R.~Schiffler.
\newblock On one-point extensions of cluster-tilted algebras.
\newblock {\em J. Algebra}, 357:168--182, 2012.

\bibitem{Rin}
C.~M. Ringel.
\newblock {\em Tame algebras and integral quadratic forms}, volume 1099 of {\em
  Lecture Notes in Mathematics}.
\newblock Springer-Verlag, Berlin, 1984.

\bibitem{SS}
R.~Schiffler and K.~Serhiyenko.
\newblock Induced and coinduced modules in cluster-tilted algebras.
\newblock {\em arXiv}, math.RT(1410.1732v2), 2014.

\bibitem{Skocrit}
A.~Skowro{\'n}ski.
\newblock Generalized standard {A}uslander-{R}eiten components without oriented
  cycles.
\newblock {\em Osaka J. Math.}, 30(3):515--527, 1993.

\end{thebibliography}
\end{document}